\documentclass[11pt]{article} 

\usepackage[utf8]{inputenc} 

\usepackage{geometry} 
\usepackage[parfill]{parskip} 

\usepackage{graphicx} 
\usepackage{array} 

\usepackage[round]{natbib}
\usepackage{hyperref,doi}
\usepackage{amsmath,amsthm,amssymb,hyperref,mathtools}
\usepackage{cases}
\usepackage{xcolor}
\usepackage{titlesec}
\usepackage{multirow}
\usepackage{qtree}

\usepackage{caption}
\usepackage{subcaption}

\titleformat{\subparagraph}
    {\normalfont\normalsize\itshape}{\thesubparagraph}{1em}{}
\titlespacing*{\subparagraph}{\parindent}{3.25ex plus 1ex minus .2ex}{.75ex plus .1ex}

\newtheorem{theorem}{Theorem}
\newtheorem{remark}{Remark}
\newtheorem{corollary}{Corollary}

\begin{document}

\title{A one-dimensional morphoelastic model for burn injuries: stability analysis, numerical validation and biological interpretation of (in)stability}

\author{	Ginger Egberts \and
           	Fred Vermolen \and
           	Paul van Zuijlen
}

\maketitle

\begin{abstract}
To deal with permanent deformations and residual stresses, we consider a morphoelastic model for the scar formation as the result of wound healing after a skin trauma. Next to the mechanical components such as strain and displacements, the model accounts for biological constituents such as the concentration of signaling molecules, the cellular densities of fibroblasts and myofibroblasts, and the density of collagen. Here we present stability constraints for the one-dimensional counterpart of this morphoelastic model, for both the continuous and (semi-) discrete problem. We show that the truncation error between these eigenvalues associated with the continuous and semi-discrete problem is of order $\mathcal{O}(h^2)$. Next we perform numerical validation to these constraints and provide a biological interpretation of the (in)stability. For the mechanical part of the model, the results show the components reach equilibria in a (non) monotonic way, depending on the value of the viscosity. The results show that the parameters of the chemical part of the model need to meet the stability constraint, depending on the decay rate of the signaling molecules, to avoid unrealistic results.
\end{abstract}

\section{Introduction}
\label{sec:1}
Burn wounds are a global problem and are the fifth most common cause of non-fatal childhood injuries. Figures show that the number of burn injuries was nearly 11 million worldwide in 2004, and about 180,000 people die from burns each year \cite{WHO}. 
Given that burns mainly occur at home and workplace and that particularly adult women and children are vulnerable to burns \cite{WHO}, targeting burn prevention specifically at these target groups results in lower numbers of incidents. 
Besides pain, itching, and loss of energy, mental factors and additional factors of wound healing play a role. Slow wound healing, infection, extreme pain, hypertrophic scars, and contractures remain as major challenges in burn management \cite{Wang}. 

The wound healing process comprises four partially overlapping phases that normally act upon each other quickly. The first phase, haemostasis, begins almost immediately after injury and aims primarily at stopping bleeding and starting the second phase. Burn wound healing passes over haemostasis, by cause of burning and cauterization of blood vessels. Hence burn wound healing starts with the second phase of normal wound healing, called the inflammatory response, which starts in just a few hours after injury to clean the wound and protects it against infections. The growth factors that play a major role stimulate angiogenesis and collagen metabolism \cite{Enoch} and activate cells, such as granulocytes (white blood cells) that play a major role in the continuation of the wound healing cascade. 

During inflammation, the wound is cleaned and protected from bacterial infections, and the proliferative phase begins. These phases in wound healing are overlapping. The sub processes that take place during the proliferative phase are re-epithelialization, angiogenesis, fibroplasia and wound contraction. Sometimes, re-epithelialization never completes and skin grafting is necessary \cite{Young}. 
The ultimate phase, remodeling and scar maturation, can take several years. This phase brings various processes and structures into balance. This results in a scar that, on average, has 50\% strength of unwounded skin (within three months), and 80\% on the long-term \cite{Enoch, Young}. 

Wound contraction is yet visible in small wounds: the edges of the wound pull in, the wound size reduces and the wounded area deforms. In adult patients, wounds can become 20-30\% smaller over several weeks \cite{Olsen1995}. Wound contraction involves a biomechanical interaction of fibroblasts, myofibroblasts, chemokines, and collagen. Depending on the wound dimensions (location on the body, size), and the extent of contraction, the result can cause reduced mobility. If the contraction result causes reduced mobility, then we commonly refer to a contracture. Contraction can lead to limited range-of-motion of joints, which can lead to immobility and is an important indication for scar revision.

Various studies report on mathematical models to predict the behavior of experimental and clinical wounds and to gain insight into which elements of the wound healing response might have a substantial influence on the contraction \cite{Tranquillo, Olsen1995, Barocas, Dallon1999,McDougall,KoppenolThesis,Menon2017} to name a few. This study uses the morphoelastic model for burn wound contraction that has been developed by Koppenol in 2017 \cite{Koppenol2017a}. Morphoelasticity is based on the following principle \cite{Hall}: the total deformation is decomposed into a deformation as a result of growth or shrinkage and a deformation as a result of mechanical forces. In a mathematical context, one considers the following three coordinate systems: ${\bf X}$, ${\bf X}_e(t)$, and ${\bf x}(t)$, which, respectively, represent the initial coordinate system, the equilibrium
at time $t$ that results due to growth or shrinkage, and the current coordinate system that results due to growth or shrinkage
and mechanical deformation. Assuming sufficient regularity, the deformation gradient tensor is written by
\begin{equation}
{\bf F} = \frac{\partial {\bf x}}{\partial {\bf X}} = \frac{\partial {\bf x}}{\partial {\bf X}_e} ~ \frac{\partial {\bf X}_e}{\partial {\bf X}}  
= {\bf A} {\bf Z},
\label{frames}
\end{equation}
in which the tensor $Z$ represents the deformation gradient tensor due to growth or shrinkage, and $A$ represents the
deformation gradient due to mechanical forces \cite{Hall,Goriely2006,Rodriguez1994}.

Given Koppenol’s morphoelastic model for skin contraction \cite{Koppenol2017a}, we analyse stability around equilibria in a one-dimensional environment to study the parametric dependence of stable and unstable solutions. We use a linear stability analysis with Fourier series, where the transformations represent perturbations around equilibria. Such a stability analysis on morphoelastic models is new in the literature. We analyse the nonlinear equations as a system of equations, and we provide stability conditions. Here we distinguish between the entire continuous problem, which represents the actual solution, and the semi-discrete problem, which is the solution of a semi-discrete solution method. We show that stability of the continuous system implies stability of the semi-discrete stable system. Next to stability conditions, we call attention to the effects of system instability regarding the real-life wound contraction. We further discuss particular components of the model that we can adapt to bring the model closer to reality. The results in this article together form an entirely recent addition to the existing morphoelastic model for skin contraction.

The organization of this paper is as follows. Section \ref{sec:2} presents the mathematical model and Section \ref{sec:3} presents the stability analysis. Subsequently, Section \ref{sec:4} presents the numerical method that is used to approximate the solution and Section \ref{sec:5} presents the numerical validation of the stability constraints and a biological interpretation of (in)stability. Finally, Section \ref{sec:6} presents the conclusion and discussion.

\section{The mathematical model}
\label{sec:2}
We borrow the morphoelastic continuum hypothesis-based modeling framework from Koppenol, and present it in one-dimensional form. It is not our aim to derive the model completely and will therefore go into this less in depth than the original articles by Koppenol \cite{Koppenol2017a, Koppenol2017b, Koppenol2017}. More details about this framework can be found in the cited articles. This model considers the displacement of the dermal layer $({u})$, the displacement velocity of the dermal layer (${v}$) and the effective strain present in the dermal layer ($\varepsilon$). The effective strain is a local measure for the difference between the current configuration of the dermal layer and a hypothetical configuration of the dermal layer where the tissue is mechanically relaxed. Furthermore, four constituents are incorporated: signaling molecules ($c$), fibroblasts ($N$), myofibroblasts ($M$) and collagen ($\rho$). Here we use collagen as a collective name for the molecules, fibrils and bundles of collagen, and we use signaling molecules as a collective name for growth factors, such as transforming growth factor beta (TGF-$\beta$), platelet derived growth factor (PDGF) and connected tissue growth factor (CTGF), and cytokines.

We show the conservation laws for mass and linear momentum, together with the evolution equation that describes how the infinitesimal effective strain changes. 
We bear in mind that due to the forces that are exerted by the cells, the domain deforms and hence the points within the domain of computation are subject to displacement. The local displacement rate is incorporated by {\em passive convection}, which is reflected by the second term in the left-hand side in equations \eqref{PDE_c}-\eqref{PDE_e}. Next, we briefly discuss what (the right-hand side of) the equations represent. 

The equation for the signaling molecules \eqref{PDE_c} represents diffusion according to normal Fickian diffusion and random spread, enhanced secretion by fibroblasts and a portion of myofibroblasts \cite{Barrientos2008}, proteolytic breakdown by Matrix Metallo Proteins (MMPs) \cite{Mast1996,Sternlicht2001}, the handle of release of MMPs by (myo)fibroblasts and collagen \cite{Lindner2012}, and the inhibition of the secretion of MMPs by signaling molecules \cite{Overall1991}:
\begin{equation}
\frac{\partial c}{\partial t} + \frac{\partial (c v)}{\partial x} = D_c \frac{\partial^2 c}{\partial x^2} + k_c \left[ \frac{c}{a_c^{I} + c} \right][N + \eta^I M] - \delta_c \frac{[N + \eta^{II}M]\rho}{1 + a_c^{II}c} c. \label{PDE_c}
\end{equation}
Here $D_c$ is the Fickian diffusion coefficient of the signaling molecules, $k_c$ is the maximum net secretion rate of the signaling molecules, $\eta^I$ is the ratio of myofibroblasts to fibroblasts in the maximum secretion rate of the signaling molecules, $a_c^{I}$ is the concentration of the signaling molecules that causes the half-maximum net secretion rate of the signaling molecules, $\delta_c$ is the proteolytic breakdown rate parameter of the signaling molecules, $\eta^{II}$ is the ratio of myofibroblasts to fibroblasts in the secretion rate of the MMPs and $1/[1+a_c^{II}c]$ represents the inhibition of the secretion of the MMPs. 
Next to the derivation of this equation in \cite{Koppenol2017a}, one finds the derivation of the second part on the right-hand side in \cite{Olsen1995}. 

The equations for the (myo)fibroblasts \eqref{PDE_N}\&\eqref{PDE_M} represent migration towards the gradient of the signaling molecules \cite{Postlethwaite1987,Boon2016,Dallon2001} by a minimal model for chemotaxis \cite{Hillen2008}, and cell density-dependent Fickian diffusion. The proliferation of the cells depends on the signaling molecules (as an activator-inhibitor), and inhibition because of crowding \cite{VandeBerg1989}. This is modeled by two similar logistic growth models. Further, the equations represent differentiation of fibroblasts to myofibroblasts \cite{Tomasek2002}, and apoptosis of the cells:
\begin{multline}
\frac{\partial N}{\partial t} + \frac{\partial (N v)}{\partial x} = -\frac{\partial}{\partial x}\left(- D_F (N+M) \frac{\partial N}{\partial x} + \chi_F N \frac{\partial c}{\partial x}\right) + \\r_F \left[ 1+\frac{r_F^{\text{max}}c}{a_c^{III}+c} \right][1-\kappa_F (N+M)] N^{1+q} - k_F c N - \delta_N N, \label{PDE_N}
\end{multline}
\begin{multline}
\frac{\partial M}{\partial t} + \frac{\partial (M v)}{\partial x} = -\frac{\partial}{\partial x}\left(- D_F (N+M)\frac{\partial M}{\partial x} + \chi_F M \frac{\partial c}{\partial x}\right) + \\r_F \left[ \frac{[1+r_F^{\text{max}}]c}{a_c^{III}+c} \right][1-\kappa_F (N+M)] M^{1+q} + k_F c N - \delta_M M. \label{PDE_M}
\end{multline}
Here $D_F$ represents (myo)fibroblast random diffusion and $\chi_F$ is the chemotactic parameter that depends on both the binding and unbinding rate of the signaling molecules with its receptor, and the concentration of this receptor on the cell surface of the (myo)fibroblasts, $r_F$ is the cell division rate, $r_F^\text{max}$ is the maximum factor of cell division rate enhancement because of the presence of the signaling molecules, $a_c^{III}$ is the concentration of the signaling molecules that cause half-maximum enhancement of the cell division rate, $\kappa_F(N+M)$ represents the reduction in the cell division rate because of crowding, $q$ is a fixed constant, $k_F$ is the signaling molecule-dependent cell differentiation rate of fibroblasts into myofibroblasts, $\delta_N$ is the apoptosis rate of fibroblasts and $\delta_M$ is the apoptosis rate of myofibroblasts. 

An important difference between the two equations is that myofibroblasts only proliferate in the presence of the signaling molecules. The form of the logistic growths needs more justification. The exact mechanism behind many of the biological processes is not always known, let alone a quantitative description of such a biological mechanism, and if others have developed a quantitative description, reliable estimates of the values for parameters are often lacking. So in general, Koppenol has avoided the use of quadratic terms in the biological parts of the models as much as possible, unless there is really a good biological reason for this. The growth of the (myo)fibroblasts is therefore taken to the power $(1+q)$ to make the model consistent. The value of $q$ is a necessary consequence of the other values of the parameters of the model. Let therefore $\overline{c},\overline{N},\overline{M}$ define the equilibria of the signaling molecules, the fibroblasts and the myofibroblasts, respectively. If we take $\overline{M}=0$ and $\overline{c}=0$ as the kinetic equilibrium, then, solving the reactive term in equation \eqref{PDE_N} for $\delta_N$ yields
\begin{equation}\label{delta_N}
\delta_N = r_F [1-\kappa_F \overline{N}] \overline{N}^{q}.
\end{equation}

The equation for collagen \eqref{PDE_R} represents the production of collagen by (myo)fibroblasts \cite{Baum2006}, enhancement of the secretion by signaling molecules \cite{Ivanoff2005}, and proteolytic breakdown of collagen by MMPs (similar as for the signaling molecules):
\begin{equation}
\frac{\partial \rho}{\partial t} + \frac{\partial (\rho v)}{\partial x} = k_\rho \left[ 1 + \left[ \frac{k_\rho^{\text{max}}c}{a_c^{IV} + c} \right] \right] [N + \eta^I M]
-\delta_\rho \frac{[N + \eta^{II}M]\rho}{1 + a_c^{II}c} \rho. \label{PDE_R}
\end{equation}
Here $k_\rho$ is the collagen secretion rate, $k_\rho^{\text{max}}$ is the maximum factor of secretion rate enhancement because of the presence of the signaling molecules, $a_c^{IV}$ is the concentration of the signaling molecules that cause the half-maximum enhancement of the secretion rate of collagen and $\delta_\rho$ is the degradation rate of collagen. A generic MMP affects the reaction kinetics of the signaling molecules and collagen, and is assumed always to be at a local equilibrium concentration. Reasoning for this modeling choice has been to avoid even more complexity and additional unknown parameter values. 

Let $\overline{\rho}$ define the equilibrium of collagen. Then, solving the reactive term in equation \eqref{PDE_R} for $\overline{\rho}$ yields
\begin{equation}\label{rho_eq}
\overline{\rho} = \sqrt{k_\rho/\delta_\rho}.
\end{equation}

The equation for the displacement velocity \eqref{PDE_v} represents Cauchy stress by a visco-elastic constitutive relation, and a body force that is proportional to the product of the cell density of the myofibroblasts and a function of the concentration of collagen. This visco-elastic constitutive relation follows the assumption from Ramtani \cite{Ramtani1,Ramtani2}, which incorporates the dependence of the Young's modulus of skin on the density of collagen:
\begin{equation}
\rho_t \left( \frac{\partial v}{\partial t} + 2v\frac{\partial v}{\partial x} \right) = \frac{\partial}{\partial x}\left( \mu\frac{\partial v}{\partial x} + E \sqrt{\rho} \varepsilon\right) + \frac{\partial}{\partial x}\left( \frac{\xi M\rho}{R^2+\rho^2} \right). \label{PDE_v}
\end{equation}
Here $\rho_t$ represents the total mass density of the dermal tissues, $\mu$ is the viscosity, $E\sqrt{\rho}$ represents the Young's modulus (stiffness), $\xi$ is the generated stress per unit cell density and the inverse of the unit collagen concentration, $R$ is a constant. The above equation represents the balance of momentum, and despite many studies neglect inertial effects (the first two terms), we have chosen to keep the inertia terms in order to stay closer to the underlying physics.
 
To incorporate a plastic deformation in the equation for the effective strain \eqref{PDE_e}, a tensor-based approach is used that is also commonly used in the context of growth of tissues (such as tumors). The ‘growth’ contribution, which with a negative sign models contraction of the tissue, is assumed to be proportional to the product of the amount of effective strain (see \cite{Hall}), the cell density of (myo)fibroblasts, and to be a function of the collagen density. In particular, it is assumed that the tensor for contraction depends on the product of the concentration of the MMPs, the concentration of the chemokines and the reciprocal of the collagen density. Taken together, the following equations present the one-dimensional morphoelastic framework for skin contraction:
\begin{equation}
\frac{\partial\varepsilon}{\partial t} + v\frac{\partial\varepsilon}{\partial x} + (\varepsilon-1)\frac{\partial v}{\partial x} = -\zeta\frac{[N+\eta^{II}M]c}{1+a_c^{II}c}\varepsilon.  \label{PDE_e}
\end{equation}
Here $\zeta$ is the rate of morphoelastic change (i.e., the rate at which the effective strain changes actively over time).

\subsection{Initial and boundary conditions}
\label{2.3}
We define the domain of computation by $\Omega_{x,t}$ and the boundary by $\partial\Omega_{x,t}$.
The dimension $x$ is in centimeters and $t$ in days. Since we are interested in the stability of the model around equilibria, we define the initial conditions by perturbations around equilibria, where the values on the boundaries are the equilibrium values. Further, we impose the following boundary conditions. For all $x\in\partial\Omega_{x,t}$ and $t\geq 0$:
\begin{equation}
\begin{split}
c(x,t) = 0, \qquad
N(x,t)=\overline{N},\qquad
M(x,t) = 0,\qquad
v(x,t) =0.
\end{split}
\end{equation}
Regarding the equations for $\varepsilon$ and $\rho$, an ordinary differential equation with derivatives regarding time in terms of the material derivative is obtained. We see this if we write the left-hand side of equation \eqref{PDE_R} as $\frac{D \rho}{D t} + \rho\frac{\partial v}{\partial x}$ and equation \eqref{PDE_e} as $\frac{D\varepsilon}{Dt}+\varepsilon\frac{\partial v}{\partial x}=\frac{\partial v}{\partial x}-\alpha\varepsilon$. The partial derivatives regarding space only involve the displacement velocity $v$. On the boundaries, for $v$ we use the boundary condition $v = 0$. Therefore, to specify the solution of $\varepsilon$ and $\rho$ in the (open) domain $\Omega$, it is unnecessary to specify any boundary conditions (the characteristics in the $x,t$-plane are vertical). We note that in cases (not currently) where characteristics would be directed out of the domain of computation, imposing these boundary conditions for would lead to failure of existence and continuity. To summarize, we do not need any boundary conditions for $\rho$ and $\varepsilon$.

\section{Linear stability of the model}
\label{sec:3}
In this Section we analyse the stability of the one-dimensional morphoelastic model for skin contraction. First, we analyse linear stability of the continuous problem. The stability conditions are formulated in terms of the input parameters. We do this analysis in order to understand the a priori behavior of the solution. Since we are not able to derive the exact solution to the problem, we also analyse stability of the numerical approximation. 
We consider the following linearised equations around equilibria $(c,N,M,\rho,v,\varepsilon)=(0,\overline{N},0,\overline{\rho},0,\overline{\varepsilon})$, where $\overline{N},\overline{\rho},\overline{\varepsilon}\in\mathbb{R}_{\geq 0}$:
\begin{equation}\label{Linearised}
\begin{split}
\frac{\partial \hat{c}}{\partial t} - D_c\frac{\partial^2\hat{c}}{\partial x^2} +\overline{N}\left[  \delta_c\overline{\rho}- \frac{ k_c }{ a_c^{I} } \right]\hat{c}=0,
\\
\frac{\partial\hat{N}}{\partial t}
\begin{multlined}[t]
- D_F\overline{N}\frac{\partial^2\hat{N}}{\partial x^2} + \chi_F\overline{N}\frac{\partial^2\hat{c}}{\partial x^2} - r_F\overline{N}^q((1+q)(1-\kappa_F\overline{N})-\kappa_F\overline{N})\hat{N} \\
+\delta_N\hat{N}+r_F\kappa_F\overline{N}^{1+q}\hat{M} - \overline{N}\left[\frac{r_Fr_F^{\text{max}}}{a_c^{III}}[1-\kappa_F\overline{N}]\overline{N}^q-k_F\right]\hat{c}=0,
\end{multlined}
\\
\frac{\partial\hat{M}}{\partial t} - D_F\overline{N}\frac{\partial^2\hat{M}}{\partial x^2} + \delta_M\hat{M} - k_F\overline{N}\hat{c} = 0,
\\
\frac{\partial \hat{\rho}}{\partial t} +\delta_\rho\overline{\rho}^2(\eta^{II}-\eta^I)\hat{M}-\delta_\rho\overline{\rho}^2\overline{N}\left(\frac{k_\rho^{max}}{a_c^{IV}}+a_c^{II}\right)\hat{c}+2\delta_\rho\overline{N}\overline{\rho}\hat{\rho}=0,
\\
\frac{\partial \hat{v}}{\partial t} - \frac{\mu}{\rho_t}\frac{\partial^2\hat{v}}{\partial x^2} - \frac{E \sqrt{\overline{\rho}}}{\rho_t}\frac{\partial\hat{\varepsilon}}{\partial x} -\frac{E\overline{\varepsilon}}{2\rho_t\sqrt{\overline{\rho}}}\frac{\partial\hat{\rho}}{\partial x} - \frac{\xi\overline{\rho}}{\rho_t(R^2+\overline{\rho}^2)}\frac{\partial\hat{M}}{\partial x}=0,
\\
\frac{\partial\hat{\varepsilon}}{\partial t} + (\overline{\varepsilon}-1)\frac{\partial\hat{v}}{\partial x} + \zeta\overline{\varepsilon}\overline{N}\hat{c}=0,
\end{split}
\end{equation}
where $\hat{c},\hat{N},\hat{M},\hat{\rho},\hat{v}$, and $\hat{\varepsilon}$ are variations around the equilibria. Here we used that $k_\rho=\delta_\rho\overline{\rho}^2$ must hold in equilibrium.

\subsection{Stability of the continuous problem} \label{3.1}
We write the variations around the equilibria in terms of a complex Fourier series
\begin{equation}\label{Variations_continuous}
\begin{alignedat}{4}
&\hat{c}(x,t) &&= \frac{1}{|\Omega|}\sum_{j = -\infty}^{\infty} c^c_j(t) e^{2 i \pi j x}, \qquad 
&&\hat{N}(x,t) &&= \frac{1}{|\Omega|}\sum_{j = -\infty}^{\infty} c^N_j(t) e^{2 i \pi j x},\\
&\hat{M}(x,t) &&= \frac{1}{|\Omega|}\sum_{j = -\infty}^{\infty} c^M_j(t) e^{2 i \pi j x}, 
&&\hat{\rho}(x,t) &&= \frac{1}{|\Omega|}\sum_{j = -\infty}^{\infty} c^{\rho}_j(t) e^{2 i \pi j x},\\
&\hat{v}(x,t) &&= \frac{1}{|\Omega|}\sum_{j = -\infty}^{\infty} c^v_j(t) e^{2 i \pi j x}, 
&&\hat{\varepsilon}(x,t) &&= \frac{1}{|\Omega|}\sum_{j = -\infty}^{\infty} c^{\varepsilon}_j(t) e^{2 i \pi j x},
\end{alignedat}
\end{equation}
where $|\Omega|$ denotes the length of $\Omega$ and $i$ represents the imaginary unit number. 

Substitution of the variations \eqref{Variations_continuous} into the linearised equations \eqref{Linearised}, multiplication by $e^{-2i\pi k x}$, and integration over $\Omega$ gives
\begin{equation}\label{Integration_continuous1}
\begin{split}
\dot{c}^c_k(t) + D_c(2\pi k)^2 c^c_k(t) +\overline{N}\left[ \delta_c\overline{\rho}- \frac{ k_c }{ a_c^{I} } \right] c^c_k(t)=0,\\
\dot{c}^N_k(t) 
\begin{multlined}[t]
+ D_F\overline{N} (2\pi k)^2 c^N_k(t) - \chi_F\overline{N} (2\pi k)^2 c^c_k(t) +r_F\kappa_F\overline{N}^{1+q} c^M_k(t)\\
- r_F\overline{N}^q((1+q)(1-\kappa_F\overline{N})-\kappa_F\overline{N}) c^N_k(t) +\delta_N c^N_k(t)\\
 - \overline{N}\left[\frac{r_Fr_F^{\text{max}}}{a_c^{III}}[1-\kappa_F\overline{N}]\overline{N}^q-k_F\right] c^c_k(t)=0,
\end{multlined}
\\
\dot{c}^M_k(t) + D_F\overline{N} (2\pi k)^2 c^M_k(t) + \delta_M c^M_k(t) - k_F\overline{N} c^c_k(t)= 0,\\
\begin{multlined}[t]
\dot{c}^{\rho}_k(t) 
+\delta_\rho\overline{\rho}^2(\eta^{II}-\eta^I) c^M_k(t)
-\delta_\rho\overline{\rho}^2\overline{N}\left[\frac{k_\rho^{max}}{a_c^{IV}}+a_c^{II}\right] c^c_k(t)\\
+2\delta_\rho\overline{N}\overline{\rho} c^{\rho}_k(t)=0,
\end{multlined}
\end{split}
\end{equation}
for the chemical part of the model, and
\begin{equation}\label{Integration_continuous2}
\begin{split}
\begin{multlined}[t]
\dot{c}^v_k(t)
+ \frac{\mu}{\rho_t} (2\pi k)^2 c^v_k(t)
- i\frac{E \sqrt{\overline{\rho}}}{\rho_t} (2\pi k) c^{\varepsilon}_k(t)
- i\frac{E\overline{\varepsilon}}{2\rho_t\sqrt{\overline{\rho}}} (2\pi k) c^{\rho}_k(t)\\
- i\frac{\xi\overline{\rho}}{\rho_t(R^2+\overline{\rho}^2)} (2\pi k) c^M_k(t)=0,
\end{multlined}
\\
\dot{c}^{\varepsilon}_k(t) + i(\overline{\varepsilon}-1) (2\pi k) c^v_k(t) + \zeta\overline{\varepsilon}\overline{N} c^c_k(t)=0,
\end{split}
\end{equation}
for the mechanical part of the model. 
The derivation of equations \eqref{Integration_continuous1} and \eqref{Integration_continuous2} is given in Appendix 1. Interchanging the second and third equation of \eqref{Integration_continuous1}, these equations together with equations \eqref{Integration_continuous2} are in the form $y' +A y = 0$ with
\begin{equation}\label{Matrix_A}
A = \begin{bmatrix}
A_{11}	&0		&0		&0		&0&0\\
A_{21}	&A_{22}	&0		&0		&0&0\\
A_{31}	&A_{32}	&A_{33}	&0		&0&0\\
A_{41}	&A_{42}	&0		&A_{44}	&0&0\\
0		&A_{52}	&0		&A_{54}	&A_{55}&A_{56}\\
A_{61}	&0		&0		&0		&A_{65}&0
\end{bmatrix}.
\end{equation}
We determine the eigenvalues of $A$ by solving $|A-\lambda I|=0$ for $\lambda$, where $I$ represents the identity matrix. For this we use the first four diagonal values as pivots and end up with a 2-by-2 matrix containing the mechanical part of the model with determinant $\lambda^2 -A_{55}\lambda -A_{56}A_{65}$. Hence the eigenvalues are the first four diagonal entries and $\lambda=\frac12 A_{55}\pm\frac12 \sqrt{A_{55}^2+4A_{56}A_{65}}$. Note that the system is linearly stable if and only if the real part of the eigenvalues is non-negative, hence we need
\begin{equation}\label{Criteria_continuous}
\begin{split}
D_c(2\pi k)^2+\overline{N}\left[  \delta_c\overline{\rho}- \frac{ k_c }{ a_c^{I} } \right]\geq0,\\
D_F\overline{N}(2\pi k)^2- r_F\overline{N}^q((1+q)(1-\kappa_F\overline{N})-\kappa_F\overline{N})+\delta_N\geq0,\\
D_F\overline{N}(2\pi k)^2 +\delta_M\geq0,\\
2\delta_\rho\overline{N}\overline{\rho}\geq0,\\
\frac{(2\pi k)^2\mu}{2\rho_t}\pm\frac12\sqrt{\left(\frac{(2\pi k)^2\mu}{\rho_t}\right)^2+4\frac{(2\pi k)^2E \sqrt{\overline{\rho}}}{\rho_t}(\overline{\varepsilon}-1)}\geq0.
\end{split}
\end{equation}
The first two requirements imply that stability is obtained for $\delta_c\geq\frac{k_c}{a_c^{I}\overline{\rho}}$ and combining the second requirement with equation \eqref{delta_N}, gives $q\delta_N\leq \kappa_Fr_F\overline{N}^{1+q}$ ($k=0$). In addition, given the relation in \eqref{delta_N}, it must hold that $\delta_N>0$ and hence $\kappa_F\overline{N}<1$. Further, the third and fourth eigenvalues meet the stability condition Re$(\lambda(A))\geq0$ independent of the chosen values for the parameters given that the parameters are positive. Lastly, linear stability is obtained for $\overline{\varepsilon} \le 1$, else a saddle point problem is obtained if $\lambda_{5,6}\in\mathbb{R}$. Note that this is also a physical requirement given that equation \eqref{PDE_e} only holds for small strains. 
These last two eigenvalues are real-valued as long as\\ $\mu \ge \frac{ \sqrt{\rho_t E\sqrt{\overline{\rho}} (1-\overline{\varepsilon})}}{\pi}$ ($k = 1$). If the last-mentioned condition is satisfied for $k = 1$, then the eigenvalues are real-valued for other values of $k$. For all the other conditions, it is obvious that they hold for all $k\in\mathbb{Z}$ as well. The constant case $k = 0$ implies that these eigenvalues are zero, which reflects the trivial case in which there are no dynamics.
This also implies that $\overline{\varepsilon} = 0$ is a stable equilibrium state with real-valued eigenvalues. 
We summarize these results in Theorem \ref{Theorem1}.

\begin{theorem} \label{Theorem1} 
Let $\{c,N,M,\rho,v,\varepsilon\}$ satisfy equations \eqref{PDE_c}-\eqref{PDE_e}. Let $\delta_N = r_F (1-\kappa_F \overline{N}) \overline{N}^{q}>0$ and $\overline{\rho} = \sqrt{k_\rho/\delta_\rho}$, then
\begin{enumerate}
\item The equilibria $(c,N,M,\rho,v,\varepsilon) = (0,\overline{N},0,\overline{\rho},0,\overline{\varepsilon})$, $\{\overline{N},\overline{\rho},\overline{\varepsilon}\}\in\mathbb{R}_{>0}$, are linearly stable if and only if $\delta_c\overline{\rho}\geq \frac{k_c}{a_c^{I}}$, and $q\delta_N\leq \kappa_Fr_F\overline{N}^{1+q}$ and $\overline{\epsilon} \le 1$;
\item Given $\overline{\varepsilon} < 1$, then 
the eigenvalues are real-valued if and only if \\
 $\mu \ge \frac{ \sqrt{\rho_t E \sqrt{\overline{\rho}} (1-\overline{\varepsilon})}}{\pi}$ ($k = 1$); 
\end{enumerate}
\end{theorem}

\begin{remark}
Note that $\delta_c \geq \frac{k_c}{a_c^{I}\overline{\rho}}$, for $k=0$ (constant states). Hence, if constant perturbations are stable, then wave-like perturbations are stable. In case $\delta_c$ is not large enough, fast oscillating perturbations will vanish, while slow oscillating perturbations will not vanish and can amplify. Further, if $\overline{\varepsilon} < 1$ and if $\mu < \frac{ \sqrt{\rho_t E \sqrt{\overline{\rho}} (1-\overline{\varepsilon})}}{\pi}$, then convergence from variations around $\overline{\varepsilon}$ will occur in a non-monotonic way over time because the eigenvalues of the linearised dynamical system are not real-valued.
\end{remark}

Next we provide some quantitative examples that illustrate the stability claims. Stability is warranted if there is a sufficient decay of the growth factor. Monotonicity (of convergence) is obtained if there is sufficient damping in terms of viscous forces.

\paragraph{Example} If we let $\delta_c=5\times10^{-4}\text{ cm}^6\text{/(cells g day)}$, $k_c=4\times10^{-13}\text{ g/(cells day)}$, $a_c^{I}=10^{-8}\text{ g/cm}^3$, and $\overline{\rho}=0.1125\text{ g/cm}^3$, then we have \\
$\delta_c=5\times10^{-4}\geq 3.55\times10^{-4}=k_c/(a_c^{I}\overline{\rho})$.
Hence with these parameter values we meet the stability condition for the signaling molecules. Further, if we let $\overline{N}=10^4<10^{6}=\kappa_F^{-1}$ cells/cm$^3$, $\delta_N=0.002$/day, $r_F=0.924$ cm$^{3q}$/(cells$^q$ day) and $q=\frac{\log(\delta_N)-\log(r_F(1-\kappa_F\overline{N})}{\log(\overline{N})}\approx-0.42$, then we have\\ $q\delta_N=-8.4\times10^{-4}\leq 1.9\times10^{-4}=\kappa_Fr_F\overline{N}^{1+q}$. Hence with these parameter values we meet the stability condition for the fibroblasts. Note that there is only a distance of $1.45\times10^{-4}\text{ cm}^6\text{/(cells g day)}$ between the left- and right-hand side in the first condition, and a much larger distance of $1.03\times10^{-3}$ between the left- and right-hand side in the second condition. In addition, substitution of $\delta_N = r_F (1-\kappa_F \overline{N}) \overline{N}^{q}$ into the second equation of \eqref{Criteria_continuous}, and solving for $q$ with $k=0$ yields $q\leq\kappa_F\overline{N}/(1-\kappa_F\overline{N})\approx0.01$, yielding the upper bound $\delta_N< 1.004$ (with the chosen parameter values). Given that the doubling time (DT) of fibroblasts ranges from 18-20 h \cite{Alberts1989,Ghosh2007}, and that the average lifespan of fibroblasts varies between 40 and 70 population doublings (PD) \cite{Ghosh2007,Moulin2011}, using the formula $\delta_N=(\ln 2)/(\text{PD}\times\text{DT}/24)$, yields the save range $0.0119\leq\delta_N\leq0.0231$ for the fibroblast apoptosis rate.

\subsection{Stability of the discrete problem} \label{3.2}
Stability of the continuous problem does not always automatically imply stability of the (semi-) discrete counterpart of the problem. Therefore, we assess stability of the semi-discrete problem, which can assess stability of the full discrete system. Lax’ Equivalence Theorem states that a consistent, stable method converges. The global truncation error tends to zero as the step size tends to zero (as $h\to0$), if the local truncation error (i.e., the difference between the derivatives and difference ratios) tends to zero as the step size is sent to zero.

A well-known way to assess numerical stability is by including Gershgorin’s Circle Theorem. This theorem is widely used and very general in the sense that it is straightforward to generalize stability to general, non-equidistant meshes and to cases where the input variables are non constant. However, in many examples, the eigenvalue bounds obtained through Gershgorin’s Circle Theorem are less accurate than by the use of the Von Neumann analysis, which is based on discrete Fourier analysis. Because of the accuracy and also the ease of application of the Von Neumann analysis, we apply this analysis on a uniform grid on the system of linearised equations with constant coefficients \eqref{Linearised}. The Von Neumann stability analysis provides sufficient conditions for numerical stability \cite{Fletcher1998}.
The \emph{finite difference method} (FDM) gives: 
\begin{equation}\label{lambda1}
\begin{split}
\lambda c_k &= - D_c\frac{c_{k-1}-2c_k+c_{k+1}}{h^2} +\overline{N}\left[  \delta_c\overline{\rho}- \frac{ k_c }{ a_c^{I} } \right]c_k,\\
\lambda N_k &\begin{multlined}[t]=  - D_F\overline{N}\frac{N_{k-1}-2N_k+N_{k+1}}{h^2}+ \chi_F\overline{N}\frac{c_{k-1}-2c_k+c_{k+1}}{h^2}\\
 +\left[\delta_N - r_F\overline{N}^q((1+q)(1-\kappa_F\overline{N})-\kappa_F\overline{N})\right]N_k+r_F\kappa_F\overline{N}^{1+q}M_k\\
- \overline{N}\left[\frac{r_Fr_F^{\text{max}}}{a_c^{III}}[1-\kappa_F\overline{N}]\overline{N}^q-k_F\right]c_k, \end{multlined}\\
\lambda M_k &= - D_F\overline{N}\frac{M_{k-1}-2M_k + M_{k+1}}{h^2} + \delta_MM_k - k_F\overline{N}c_k,
\\
\lambda \rho_k &=\delta_\rho\overline{\rho}^2(\eta^{II}-\eta^I)M_k - \delta_\rho\overline{\rho}^2\overline{N}\left(\frac{k_\rho^{max}}{a_c^{IV}}+a_c^{II}\right)c_k + 2\delta_\rho\overline{N}\overline{\rho}\rho_k,
\end{split}
\end{equation}
for the chemical part of the model, and
\begin{equation}\label{lambda2}
\begin{split}
\lambda v_k &\begin{multlined}[t]= -\frac{\mu}{\rho_t}\frac{v_{k-1} - 2v_k + v_{k+1}}{h^2} - \frac{E\sqrt{\overline{\rho}}}{\rho_t}\frac{\varepsilon_{k+1}-\varepsilon_{k-1}}{2h}\\
- \frac{E\overline{\varepsilon}}{2\rho_t\sqrt{\overline{\rho}}}\frac{\rho_{k+1}-\rho_{k-1}}{2h}
-\frac{\xi\overline{\rho}}{\rho_t(R^2+\overline{\rho}^2)}\frac{M_{k+1}-M_{k-1}}{2h},\end{multlined}\\
\lambda\varepsilon_k &= (\overline{\varepsilon}-1)\frac{v_{k+1}-v_{k-1}}{2h}+\zeta\overline{\varepsilon}\overline{N}\hat{c}_k,
\end{split}
\end{equation}
for the mechanical part of the model. 
Let
\begin{equation}\label{Variations_discrete}
\begin{alignedat}{6}
  &c_k 		&&=\sum_{\beta=1}^{n-1} \hat{c}_\beta e^{-2\pi\beta khi}, \quad &&N_k &&=\sum_{\beta=1}^{n-1} \hat{N}_\beta e^{-2\pi\beta khi}, \quad &&M_k &&=\sum_{\beta=1}^{n-1} \hat{M}_\beta e^{-2\pi\beta khi},\\
  &\rho_k 	&&=\sum_{\beta=1}^{n-1} \hat{\rho}_\beta e^{-2\pi\beta khi}, &&v_k &&=\sum_{\beta=1}^{n-1} \hat{v}_\beta e^{-2\pi\beta khi}, &&\varepsilon_k &&=\sum_{\beta=1}^{n-1} \hat{\varepsilon}_\beta e^{-2\pi\beta khi}.
\end{alignedat}
\end{equation}
Substitution of \eqref{Variations_discrete} in equations \eqref{lambda1} and \eqref{lambda2}, subdivision by $e^{-2\pi\beta khi}$, and using Euler's formula and $2-2\cos(2\pi\beta h)=4\sin^2(\pi\beta h)$ results in
\begin{equation}\label{Substitution_discrete1}
\begin{split}
\lambda \hat{c}_\beta & 
\begin{multlined}[t]=
\frac{D_c}{h^2}4\sin^2(\pi\beta h)\hat{c}_\beta +\overline{N}\left[\delta_c\overline{\rho} - \frac{k_c}{a_c^{I}}\right]\hat{c}_\beta,
\end{multlined}\\
\lambda \hat{N}_\beta &
\begin{multlined}[t]=  
\frac{D_F\overline{N}}{h^2}4\sin^2(\pi\beta h)\hat{N}_\beta
- \frac{\chi_F\overline{N}}{h^2}4\sin^2(\pi\beta h)\hat{c}_\beta\\
+\left[\delta_N - r_F\overline{N}^q((1+q)(1-\kappa_F\overline{N})-\kappa_F\overline{N})\right]\hat{N}_\beta\\
+r_F\kappa_F\overline{N}^{1+q}\hat{M}_\beta
- \overline{N}\left[\frac{r_Fr_F^{\text{max}}}{a_c^{III}}[1-\kappa_F\overline{N}]\overline{N}^q-k_F\right]\hat{c}_\beta, 
\end{multlined}\\
\lambda \hat{M}_\beta &
\begin{multlined}[t]= 
\left[\frac{D_F\overline{N}}{h^2}4\sin^2(\pi\beta h) +\delta_M\right]\hat{M}_\beta - k_F\overline{N}\hat{c}_\beta,
\end{multlined}\\
\lambda \hat{\rho}_\beta &
\begin{multlined}[t]=
\delta_\rho\overline{\rho}^2(\eta^{II}-\eta^I)\hat{M}_\beta
- \delta_\rho\overline{\rho}^2\overline{N}\left[\frac{k_\rho^{max}}{a_c^{IV}}+a_c^{II}\right]\hat{c}_\beta
+ 2\delta_\rho\overline{N}\overline{\rho}\hat{\rho}_\beta,
\end{multlined}
\end{split}
\end{equation}
for the chemical part of the model, and
\begin{equation}\label{Substitution_discrete2}
\begin{split}
\lambda \hat{v}_\beta &
\begin{multlined}[t]= 
\frac{\mu}{\rho_t h^2}4\sin^2(\pi\beta h)\hat{v}_\beta
+i \frac{E\sqrt{\overline{\rho}}}{\rho_t h}\sin(2\pi\beta h)\hat{\varepsilon}_\beta\\
+i \frac{E\overline{\varepsilon}}{2\rho_t\sqrt{\overline{\rho}}h}\sin(2\pi\beta h)\hat{\rho}_\beta
+i \frac{\xi\overline{\rho}}{2\rho_t(R^2+\overline{\rho}^2)h}\sin(2\pi\beta h)\hat{M}_\beta,\end{multlined}\\
\lambda\hat{\varepsilon}_\beta &
\begin{multlined}[t]= 
-i\frac{(\overline{\varepsilon}-1)}{h}\sin(2\pi\beta h)\hat{v}_\beta
+\zeta\overline{\varepsilon}\overline{N}\hat{c}_\beta,\end{multlined}
\end{split}
\end{equation}
for the mechanical part of the model. The derivation of equations \eqref{Substitution_discrete1} and \eqref{Substitution_discrete2} is given in Appendix 2. These equations are in the form $\lambda z=C z$ with the matrix $C$ as in \eqref{Matrix_A}. Hence the eigenvalues are found in the same way as before. Note that the discrete system is linearly stable if and only if the real part of the eigenvalues is non-negative, hence we need
\begin{equation}\label{Criteria_discrete}
\begin{split}
\frac{D_c}{h^2}4\sin^2(\pi\beta h) +\overline{N}\left[\delta_c\overline{\rho} - \frac{k_c}{a_c^{I}}\right]\geq0,\\
\frac{D_F\overline{N}}{h^2}4\sin^2(\pi\beta h)- r_F\overline{N}^q((1+q)(1-\kappa_F\overline{N})-\kappa_F\overline{N})+\delta_N \geq0,\\
\frac{D_F\overline{N}}{h^2}4\sin^2(\pi\beta h) +\delta_M\geq0,\\
2\delta_\rho\overline{N}\overline{\rho}\geq0,\\
\begin{multlined}[t]
\frac{2\mu}{\rho_th^2}\sin^2(\pi\beta h)\pm
\\\frac12\sqrt{\left(\frac{\mu}{\rho_t h^2}4\sin^2(\pi\beta h)\right)^2+4\frac{E\sqrt{\overline{\rho}}}{\rho_t h^2}(\overline{\varepsilon}-1)\sin^2(2\pi\beta h)}\geq0.
\end{multlined}
\end{split}
\end{equation}
To guarantee linear stability, the first requirement states $\delta_c\overline{\rho}\geq \frac{k_c}{a_c^{I}}$. Given $\delta_N = r_F (1-\kappa_F \overline{N}) \overline{N}^{q}$, the second requirement states $q\delta_N\leq \kappa_Fr_F\overline{N}^{1+q}$. The third and fourth eigenvalues meet the stability condition independent of the chosen values for the parameters given that the parameters are positive. Lastly, for the discrete problem linear stability is also obtained for $\overline{\varepsilon} \le 1$, and since
\begin{equation}
4\frac{E\sqrt{\overline{\rho}}}{\rho_t h^2}(\overline{1-\varepsilon})\sin^2(2\pi\beta h)\geq0,
\end{equation}
stability is guaranteed for all $h\in\mathbb{R}_{>0}$. To conclude, we have demonstrated that if the equilibrium is stable in the continuous problem, then it is also stable in the semi-discrete problem.

There exists a consistency between the stability criteria of the continuous problem and the stability criteria of the discrete problem. We show this consistency by writing $\sin^2(x)$ as a Taylor series. Substitution into the first and last equation in \eqref{Criteria_discrete} yields:
\begin{equation}
\begin{split}
D_c(2\pi\beta)^2 +\mathcal{O}(h^2) +\overline{N}\left[\delta_c\overline{\rho} - \frac{k_c}{a_c^{I}}\right]\geq0,\\
\begin{multlined}[t]
\frac{(2\pi\beta)^2\mu}{2\rho_t}+\mathcal{O}(h^2)\pm
\\\frac12\sqrt{\left(\frac{(2\pi\beta)^2\mu}{\rho_t}+\mathcal{O}(h^2)\right)^2+4\frac{(2\pi\beta)^2E\sqrt{\overline{\rho}}}{\rho_t}(\overline{\varepsilon}-1)+\mathcal{O}(h^2)}\geq0.\end{multlined}
\end{split}
\end{equation}
Comparison to the first and last second equation of \eqref{Criteria_continuous}
\begin{equation}
\begin{split}
D_c(2\pi k)^2+\overline{N}\left[  \delta_c\overline{\rho}- \frac{ k_c }{ a_c^{I} } \right]\geq0,\\
\frac{(2\pi k)^2\mu}{2\rho_t}\pm\frac12\sqrt{\left(\frac{(2\pi k)^2\mu}{\rho_t}\right)^2+4\frac{(2\pi k)^2E \sqrt{\overline{\rho}}}{\rho_t}(\overline{\varepsilon}-1)}\geq0.
\end{split}
\end{equation}
yields a difference in eigenvalues of order $\mathcal{O}(h^2)$. Note that in the same way, a difference of order $\mathcal{O}(h^2)$ follows for the second equations of \eqref{Criteria_continuous} and \eqref{Criteria_discrete}. 

Furthermore, the last equation in \eqref{Criteria_discrete} implies that for real-valued eigenvalues we need
\begin{equation*}
\frac{\mu^2}{\rho_t^2 h^4}4^2\sin^4(\pi\beta h) \geq 4\frac{E\sqrt{\overline{\rho}}}{\rho_t h^2}(1-\overline{\varepsilon})\sin^2(2\pi\beta h).
\end{equation*}
Writing $\sin^2(2\pi\beta h)=4\sin^2(\pi\beta h)\cos^2(\pi\beta h)$, multiplication by $\frac{\rho_t^2h^2}{4^2\sin^4(\pi\beta h)}$ gives
\begin{equation*}
\mu^2\geq \rho_t h^2E\sqrt{\overline{\rho}}(1-\overline{\varepsilon})\frac{\cos^2(\pi\beta h)}{\sin^2(\pi\beta h)}.
\end{equation*}
Hence the numerical criterium
\begin{equation}
\mu \geq \frac{h}{\tan(\pi\beta h)} \sqrt{\rho_tE\sqrt{\overline{\rho}}(1-\overline{\varepsilon})}.
\end{equation}
For consistency, we have
$$
\lim_{h\to0}\frac{h}{\tan(\pi\beta h)}=\lim_{h\to0}\frac{\pi\beta h}{\tan(\pi\beta h)}\cdot\frac{1}{\pi\beta}=\frac{1}{\pi\beta}
$$
and $\frac{h}{\tan(\pi \beta h )} \le \frac{1}{\pi \beta}$, for $\beta = 1,\ldots, n-1$ ($h n = |\Omega|$). Hence for monotonic convergence for $\beta=1$, we see that the convergence is consistent with convergence of the fully continuous model for $h\to0$. We summarise the results in Theorem \ref{Theorem2}.
\begin{theorem} \label{Theorem2}
Let $\{c,N,M,\rho,v,\varepsilon\}$ satisfy the semi-discrete spatial finite differences version of equations \eqref{PDE_c}-\eqref{PDE_e}. Then stability in the fully continuous problem implies stability in the semi-discrete formulation, regardless the step-size. Furthermore, monotonic convergence in the fully continuous problem implies monotonic convergence in the semi-discrete problem formulation, regardless the step-size.
\end{theorem}
\begin{corollary}
Let $\{c,N,M,\rho,v,\varepsilon\}$ satisfy the semi-discrete spatial finite differences version of equations \eqref{PDE_c}-\eqref{PDE_e}. Let $\delta_N = r_F (1-\kappa_F \overline{N}) \overline{N}^{q}$ and $\overline{\rho} = \sqrt{k_\rho/\delta_\rho}$, then the equilibria are unconditionally stable for the trapezoid rule and the Euler backward method as long as $\delta_c\overline{\rho}\geq {k_c}/{a_c^{I}}$ and $q\delta_N\leq \kappa_Fr_F\overline{N}^{1+q}$. Furthermore, the Euler backward method is A-stable.
\end{corollary}
\begin{remark}
It is possible that the semi-discrete yields monotonic convergence, whereas the continuous problem does not. The reason for this is that $\frac{h}{\tan(\pi\beta h)} \leq \frac{1}{\pi\beta}$. Hence the inequality for the continuous problem is sharper than for the semi-discrete problem.
\end{remark}

\section{Numerical method for validation}
\label{sec:4}
We approximate the solution to the model equations by the finite-element method using linear basis functions. For more information about this method we refer to \cite{VanKan}. We multiply the equations \eqref{PDE_c}-\eqref{PDE_e} by a test function $\varphi(x,t)\in H_0^1$, integrate over the domain of computation $\Omega$ (integration by parts), apply the application of the Gauss' theorem, and apply the Leibniz-Reynold\rq{}s transport theorem. 

To construct the basis functions we subdivide the domain of computation into $n\in\mathbb{N}$ sub-domains $e_p=[x_p,x_{p+1}]$ (i.e., the elements). Let $X_h(t) = \bigcup e_p$ the finite element subspace and $x_j , j \in \{1,\dots,n+1\}$ the vertices of the elements. We choose $\varphi_i(x_j,t)=\delta_{ij}$, $i,j\in\{1,\dots,n+1\}$ as the linear basis functions, where $\delta_{ij}$ denotes the Kronecker delta function. 

Note that the following holds for the chosen subspace $X_h(t) \subset \Omega_{x,t}$:\\ $\frac{\mathrm{D}\varphi_i}{\mathrm{D}t}=0$ for all $\varphi_i$ \cite{Dziuk}. The Galerkin equations are simplified using this property. We solve the Galerkin equations using backward Euler time integration and we use a monolithic approach with inner Picard iterations to account for the non-linearity of the equations. To avoid loss of monotonicity (i.e. oscillations), we use the process called mass lumping. 

We approximate the local displacements by
\begin{equation}
u_i^{t+\Delta t}\simeq u_t^t+\Delta t v_i^{t+\Delta t},
\end{equation}
with
\begin{equation}
u(x,0)=0,\quad\forall x\in\Omega_{x,0},
\end{equation}
the initial condition.

\section{Results}
\label{sec:5}
To experimentally assess the convergence of the numerical method, we use a domain of computation of 10 cm in which we model a 4 cm large wound. To account for the steepness of the gradients of the initial fibroblast distribution, and signaling molecule and collagen densities, we use an interval with length of 1 cm over which the initial solution varies between its equilibrium and zero. Within the wound, we assume that there are 2000 fibroblast cells/cm$^3$, $10^{-8}$ g/cm$^3$ signaling molecules and $0.01125$ g/cm$^3$ collagen present. We model the gradient of the steepness area by sine functions. We divide the domain of computation in $n$ elements, where $n\in\{41, 81, 161, 321, 641, 1281\}$. For each simulation, we define $\Delta t=h^2$, where $h$ is the size of the elements, and simulate skin contraction for 1 day. In each simulation, we report the densities of the variables (the solutions) and the relative surface area of the wound (RSAW). The convergence order results are computed as follows. 
Let $\lim\limits_{h\to0}z_h(x,1)=z(x,1)$ denote the true density of variable $z$ on day 1 and $z_{0.0078}(x,1)=:z_{h/r}$ the solution in the last simulation (i.e., the reference, which has been computed using the highest numerical resolution). We approximate the errors $\epsilon:=\int|z-z_h|\mathrm{d}x$ of the solutions on the full domain of computation, and since we are interested in the displacement of the boundary of the wound, we approximate the errors of the solutions on the boundary of the wound in particular. For this we use the following error definition:
\begin{equation}\label{norm_e41}
\epsilon_{|41|}(h) =  \sum_{i=1}^{41} \left|z_{h/r}(x_{i,41})-z_h(x_{i,41})\right|,
\end{equation}
where the grid-points $x_{i,n}$ correspond to the grid-points in the simulation with $n=41$ nodes. This error is a variant of the $L^1$-norm in which we evaluate the solution to the equations on the same grid-points. 
 
\begin{figure}[!ht]
\centering
\includegraphics[width=\textwidth]{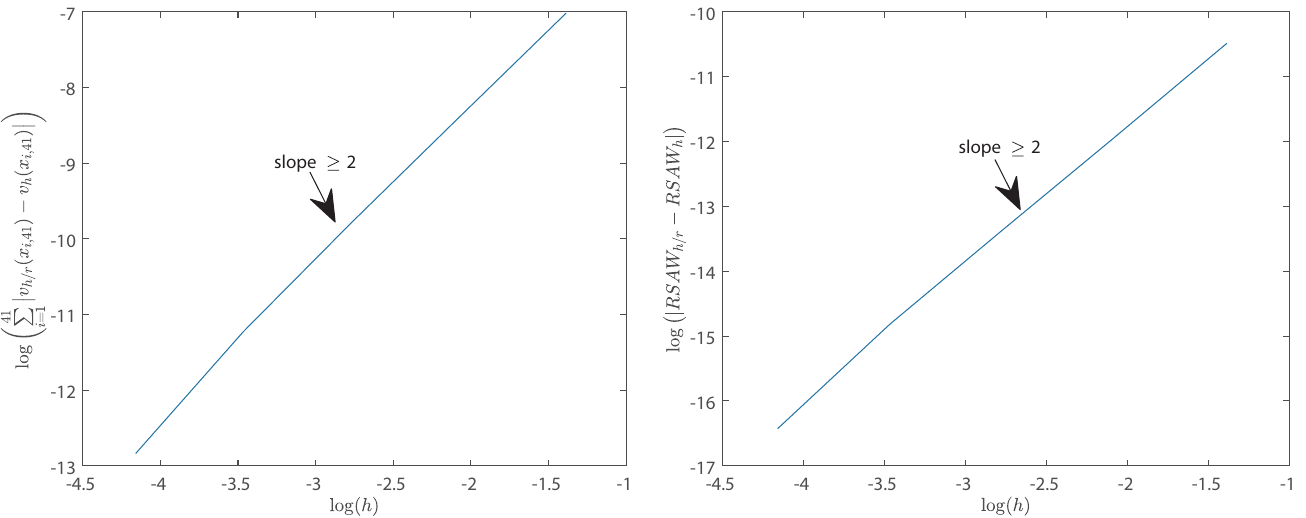}
\caption{Numerical validation of convergence. Here, the contraction of a wound of 4 cm, with 1 cm steepness, on a domain of 10 cm is simulated. Initially, in the wound there are 2000 fibroblast cells/cm$^3$, $10^{-8}$ g/cm$^3$ signaling molecules and $0.01125$ g/cm$^3$ collagen present. The values of the other parameters are shown in Table \ref{Parameters}. The left plot shows the logarithm of the step size versus the logarithm of the absolute displacement velocity density error on the domain of computation on a fixed number of grid-points. The right plot shows the logarithm of the step size $h$ versus the logarithm of the absolute relative surface area error}
\label{fig:convergence}
\end{figure}

Figure \ref{fig:convergence} shows some of the results for error $\epsilon_{|41|}$, where we show the relations of the errors with the element size $h$ for the displacement velocity, and the error of the relative surface area of the wound.

From the left plot, we see that the absolute error of the displacement velocity decreases consistently as $h$ becomes smaller. The average slope of this graph is 2.1882, hence the order of convergence is about $\mathcal{O}(h^2)$. From the right plot, we see that the absolute error of the relative surface area of the wound decreases consistently as $h$ becomes smaller. The average slope of this graph is 2.2092, showing an order of convergence about $\mathcal{O}(h^2)$ as well. We note that all averaged slopes of the logarithms of the  absolute errors of the variables and the relative surface area of the wound show an overall consistent convergence of order $\mathcal{O}(h^2)$. One finds these slopes in Table \ref{tab:slopes} in Appendix 4.

To validate the stability for the model, we perturb the initial conditions around equilibria using sine functions, and we vary the parameters $\delta_c$ and $\mu$. 
We use $n=500$ elements to divide the domain of computation between 0 and 1, which represents half a domain of the modeled skin on which we perform computations. This is possible due to the symmetry of the model. We fix all parameters except for $\delta_c$ and $\mu$. The values of the fixed parameters are given in Table \ref{Parameters}.
When not stated otherwise, for time integration, we use a step of $\Delta t=5\times10^{-1}$ days. 

For the initial conditions, we vary the number of waves $k$ using three levels (1, 5 and 10).
We perturb the initial condition for the fibroblasts and collagen by using a sine function with amplitude $10$ cells/cm$^3$ and $10^{-2}$ g/cm$^3$, respectively. 
This is possible because the equilibrium distribution of the fibroblasts and the equilibrium density of collagen are non-zero. 
For the initial condition of the myofibroblasts and the signaling molecules we use uniform splines with $2k+1$ knots. On the boundaries, the knots have zero value, and in between the values are 3 and 6 cells/cm$^3$ for the myofibroblasts, and $0.5\times10^{-15}$ and $2\times10^{-15}$ g/cm$^3$ for the signaling molecules. This way we ensure that the myofibroblast distribution and signaling molecule density values are positive. 
The initial amplitudes of the displacement velocity and effective strain are 0.05 and 0.5, respectively.

For stability, Theorem \ref{Theorem1} requires that $\delta_c\geq\frac{k_c}{a_c^{I}\overline{\rho}}$ in case $k=0$. 
Further, given that the equilibrium density of the effective strain is less than 1, eigenvalues are real-valued if and only if $\mu \ge \sqrt{\rho_t E \sqrt{\overline{\rho}} (1-\varepsilon_0)}/\pi$ in case $k = 1$. We choose to vary the signaling molecule decay rate $\delta_c$ using three levels ($2\times10^{-4}$, $3\times10^{-4}$ and $5\times10^{-4}$) cm$^6$/(cells g day), where the first two values are chosen such that such that the stability condition is not met. We vary the viscosity parameter $\mu$ using two levels (1 and 100) (N day)/cm$^2$. The first value is chosen such that the corresponding eigenvalue is not real-valued. Videos corresponding to the shown figures can be found in the \href{https://figshare.com/s/0875ae5579bbe17f352f}{online resources}. 
The values of the fixed parameters can be found in Appendix 3. 

\begin{figure}[!ht]
\centering
\includegraphics[width=\textwidth]{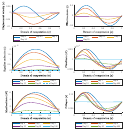}
\caption{Evolution of distributions and densities of the modeled variables for $\delta_c=5\times10^{-4}$ cm$^6$/(cells g day) and $\mu=100$ (N day)/cm$^2$. The values of the other parameters are shown in Table \ref{Parameters}. The plots on the upper left and right, the middle left and right, and the lower left and right show the displacement velocity, the effective strain, the signaling molecules, the fibroblasts, the myofibroblasts, and collagen, respectively}
\label{Stable1}
\end{figure}

In the first simulation we take $\delta_c=5\times10^{-4}$ cm$^6$/(cells g day) and $\mu=100$ (N day)/cm$^2$ and simulate for 400 days. We note that for these values, the stability criteria are met. Figure \ref{Stable1} shows the results. 

We see that the displacement velocity density rearranges to negative values. As the density moves below zero, the amplitude of the wave initially increases, after which the density moves gradually toward the equilibrium $v=0$. The effective strain density does not change signs. The values on the boundaries of the domain of computation initially move away from the equilibrium, where all other values gradually move toward the equilibrium $\varepsilon\approx-0.05$. Because of the boundary condition, the signaling molecule density is fixed at equilibrium on the left boundary of the domain of computation. We see that on the right boundary, the density increases in the first days, after which it decreases to the equilibrium $c=0$ g/cm$^3$. Because of the negative values of the displacement velocity density after 12 hours, the mesh moves to the left. This is most clear in the fibroblast plot. During the simulation, the fibroblast distribution displacements to the left, and values above the equilibrium gradually move toward the equilibrium $N=10^4$ cells/cm$^3$. The fibroblast distribution on the right boundary starts by moving away from the equilibrium as the fibroblasts differentiate to myofibroblasts because of the increased density of signaling molecules. After the signaling molecule density is almost zero around the right boundary on day 30, the fibroblast distribution moves toward the equilibrium, reaching it fully around day 400. We see the same effect in the myofibroblast plot, where the distribution moves to the left, and moves gradually toward the equilibrium $M=0$ cells/cm$^3$. Only the values on the right boundary move away from the equilibrium in the first 10 days, because of the differentiated fibroblasts. The plot of collagen is like the plot of the effective strain, although the effect of the local displacements seems larger for collagen, and for collagen it takes much longer before the density reaches the equilibrium $\rho=0.1125$ g/cm$^3$. Overall, the model behaves absolutely stable given these stable parameter values. 

From a biological perspective, minor variations in the number of (myo) fibroblast cells, and in the density of signaling molecules and collagen, already initialises wound healing in which contraction appears for 100 days. If there is a disruption in the distribution of collagen, the skin recovers this almost immediately. However, this process takes longer than for signaling molecules, for example. Further, local displacements in the skin are either in the direction toward the center of the wound or in the direction of the boundary of the wound. 

Next, in the second, third and fourth ($k=1,5,10$) simulations we take $\delta_c=2\times10^{-4}$ cm$^6$/(cells g day) and $\mu=100$ (N day)/cm$^2$ and simulate for 1200 days (not shown). While running these simulations, at first the constituents (almost) reach their equilibria. For $k=1$, the signaling molecule density reaches the equilibrium around day 250, the fibroblast distribution changes towards equilibrium until day 650, the myofibroblast distribution reaches equilibrium around day 390, and the collagen density around day 650 as well. Further, both the displacement velocity and effective strain density reach equilibria within 15 days. However, from day 660, the signaling molecule density increases and starts decreasing around day 753. The fibroblast distribution decreases after day 650 and starts increasing around day 745 again. The myofibroblast distribution also increases, which happens around day 638, and starts decreasing again around day 704. Shortly after the collagen density seems to reach equilibrium around day 650, the density explodes and does not start decreasing. 
Because of singular matrices, we ended this simulation. The Picard iterations did not converge and because of NaN’s in all the solutions, there were no plots available anymore. The same is seen where $k=5,10$. 

Theoretically, if the human body or an external factor reduces the decay rate of signaling molecules too much, then initially, this does not cause the skin to rupture. However, after a few years, the secretion of signaling molecules can increase significantly, causing such problems. Present fibroblasts fully differentiate into myofibroblasts. The scar will undergo a severe contraction, and collagen will cause tissue to rupture because of excessive production.
We believe that the human body protects against the lowering of the decay rate of signaling molecules to this extent, in order to prevent such a non-realistic occurrence.

\begin{figure}[!ht]
\centering
\includegraphics[width=\textwidth]{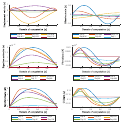}
\caption{Evolution of distributions and densities of the modeled variables for $\delta_c=5\times10^{-4}$ cm$^6$/(cells g day) and $\mu=1$ (N day)/cm$^2$. The values of the other parameters are shown in Table \ref{Parameters}. The plots on the upper left and right, the middle left and right, and the lower left and right show the displacement velocity, the effective strain, the signaling molecules, the fibroblasts, the myofibroblasts, and collagen, respectively}
\label{Stable2}
\end{figure}

Next, in the fifth, sixth and seventh ($k=1,5,10$) simulations we take $\delta_c=5\times10^{-4}$ cm$^6$/(cells g day) and $\mu=1$ (N day)/cm$^2$ and simulate for 600 days. Note that the signaling molecule decay rate stability condition is met and that we focus on the effect of complex eigenvalues in the mechanical part of the model.  Initially we use a time step of $\Delta t=0.01$, and we change that to $\Delta t=1$ after 2 days, and to $\Delta t=2$ after 50 days. Figure \ref{Stable2} shows the results for $k=1$. 

We see that all the constituents reach equilibria within 600 days, after which the distributions and densities do not change anymore. Initially, the displacement velocity density oscillates around zero, moving the mesh to the left and right, and the effective strain density oscillates around the (new) equilibrium. Shortly after the start of the simulation, the wave in the displacement velocity density fades out. Further, within approximately 15 minutes, the amplitude increases by a factor 10 above the equilibrium value, and by a factor 25 below the equilibrium value. Shortly after that, around approximately 1.5 hours, the amplitude of the displacement velocity density has increased by a factor 45, after which the amplitude decreases until zero. Both the displacement velocity density and effective strain density reach the equilibria within a few days, the displacement velocity density reaching the equilibrium $v=0$ first. Note that these results both confirm the non-monotonic convergence from the variations around $\overline{\varepsilon}$ (see Theorem \ref{Theorem1} and Theorem \ref{Theorem2}).
We see the mesh also moving in the plots of the constituents. While the displacement velocity density oscillates, the distributions and densities of the constituents move from the right to the left and back, until the distributions and densities move gradually towards the equilibria. First, the signaling molecules density reaches equilibrium around day 60. About twice that time, around 120 days, the myofibroblast distribution reaches equilibrium. The fibroblast distribution grows as follows. After a few days, when the displacement velocity density reaches equilibrium, the fibroblast distribution above the equilibrium decreases, and the fibroblast distribution below the equilibrium increases, except for the fibroblast distribution around the right boundary of the domain of computation, representing the center of the portion of skin that we model. The number of fibroblasts around this right boundary decreases until about 23 days, after which it increases towards equilibrium. The collagen density changes calmly: the density above the equilibrium moves downward to the equilibrium, and the density below the equilibrium moves upward to the equilibrium. 

Where $k=5$ (figures not shown), the results show that increasing the number of waves makes the initial increase in amplitudes in the displacement velocity density smaller. Again, initially this amplitude increases around 15 minutes, after which it decreases while the density oscillates around the equilibrium. Fading out the waves takes more time, here about 4.8 hours, compared to 1.5 hours where $k=1$, and the local displacements are much smaller. The other densities and distributions change similar to where $k=1$, except for some features. Equilibria are reached around day 112, 210, and 600 for the signaling molecule density, the myofibroblast distribution, the fibroblast distribution and collagen density, respectively, the first two later than where $k=1$. The waves in the fibroblast distribution disappear faster and the distribution moves faster toward the equilibrium. 
The smaller local displacements are clearly visible in the plots of the constituents. We have seen that the signaling molecule density shifts to the left between 0 and 3 hours, and to the right between 3 and 8 hours. Further, the density decreases gradually toward the equilibrium, and the waves have already started fading out on day 4. 

Comparing the results from the simulation where $k=10$ (figures not shown) with the simulations where $k=1,5$, we conclude the waves fade out faster for faster oscillating perturbations and that initially the distributions and densities of the constituents and the effective strain move faster toward the equilibria. In addition, the initial increase in amplitude in the displacement velocity density is larger for smaller $k$. Taken these numerical results together, we can confirm that the one-dimensional morphoelastic framework for skin contraction is stable given that $\delta_c\geq k_c/(a_c^{I}\overline{\rho})$.

From a biological perspective, a large value of the viscosity mimics a large amount of damping, and this damping term makes the equation for the displacement velocity more ‘diffusive’. A diffusion equation satisfies a maximum principle, that is, the extremes can only be assumed on the boundary of the domain or initially, unless the solution is constant. This implies that the solution must behave more monotonically for large viscosities as shown in the upper left plot in Figure \ref{Stable1}. A small value in the viscosity makes the equation for the displacement velocity less diffusive, so that the extremes are not bounded by the boundary of the domain or initially as shown in the upper left plot in Figure \ref{Stable2}. Here, the modeled medium is less resistant to the rate of deformation, and given the initial fluctuation in the displacement velocity density, this results not only in a back-and-forth movement in the displacement, also a direct effect in the stress (effective strain) that is proportional to the shear deformation.

\begin{figure}[!ht]
\centering
\includegraphics[width=\textwidth]{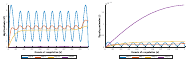}
\caption{(early) Evolution of myofibroblast distribution and signaling molecule density for $\delta_c=3\times10^{-4}$ cm$^6$/(cells g day) and $\mu=100$ (N day)/cm$^2$ and $k=10$.  The values of the other parameters are shown in Table \ref{Parameters}. The left and right plots show the myofibroblasts and the signaling molecules, respectively}
\label{Stable3}
\end{figure}

\begin{figure}[!ht]
\centering
\includegraphics[width=\textwidth]{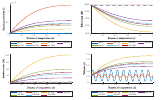}
\caption{Evolution of distributions and densities of the constituents for $\delta_c=3\times10^{-4}$ cm$^6$/(cells g day) and $\mu=100$ (N day)/cm$^2$ and $k=10$. The values of the other parameters are shown in Table \ref{Parameters}. The plots on the upper left and right, and lower left and right, show the signaling molecules, the fibroblasts, the myofibroblasts and collagen, respectively}
\label{Stable4}
\end{figure}

As stated before, the model can numerically be unstable when $\delta_c< k_c/(a_c^{I}\overline{\rho})$. However, we have seen that sometimes for small signaling molecule decay rates not too far below the stated lower boundary, the model still converges. In the last simulation we set the number of waves with $k=10$, and we take $\delta_c=3\times10^{-4}$ cm$^6$/(cells g day) and $\mu=100$ (N day)/cm$^2$. Figures \ref{Stable3} and \ref{Stable4} show some of the results of the simulation of 1000 days. These show that the model converges, and highlight what happens in this case. 

First, everything seems calm until day 60 (for example, see the left plot in Figure \ref{Stable3}). The displacement velocity density (figure not shown) reaches equilibrium within 10 days, and the effective strain density around day 20. However, the initial perturbed waves are still visible. Initially, the signaling molecule density decreases, but on approximately day 9 the upper bound of the density surpasses the initial upper bound (see the right plot in Figure \ref{Stable3}). The signaling molecule density keeps increasing until day 215, affecting the (myo)fibroblast distributions and the collagen density, shown in Figure \ref{Stable4}. 

The initial perturbed waves in the (myo)fibroblast distribution fade out within 4.5 days. Both distributions move toward the corresponding equilibria $10^4$ cells/c$^3$ and approximately 0.16 cells/cm$^3$ (hence no cells), respectively. However, on days 63.5 and 65, for the fibroblasts and myofibroblasts respectively, the distributions move away from the equilibria. Only the collagen density is not affected by this setup until day 120, after which this density increases. 

After the signaling molecule density decreases from day 215 on, the myofibroblast distribution, the collagen density and the fibroblast distribution keep moving away from their equilibria until days 230, 250 and 260, respectively. From the plot for collagen, it takes more time to fade out the initial perturbed waves. From these moments (i.e., the days where maxima and minima are reached), the distributions and densities of the constituents oscillate around a new equilibrium. At the end of the simulation of 1000 days the new equilibria in the center of the modeled skin are $4.245\times10^{-11}$ g/cm$^3$, 9723 cells/cm$^3$, 76 cells/cm$^3$, and 0.1348 g/cm$^3$ for the signaling molecules, the fibroblasts, the myofibroblasts, and collagen, respectively.

From a biological perspective, if there is an enhanced expression of signaling molecules because of their reduced decay, a wound may heal properly at first. However, over time, persistent signaling will lead to over-expression of signaling molecules, resulting in excessive scarring and contraction.
The excessive deposition of collagen is reminiscent of keloids and hypertrophic scars, characterised by thicker collagen bundles \cite{Tuan1998}. In addition, myofibroblasts are abundant in hypertrophic scars. Since aberrant TGF-$\beta$ signaling in myofibroblasts is associated with the formation of hypertrophic scars \cite{Zhang2020}, it is likely that such a situation exists precisely because of a lower decay rate of signaling molecules. Further, hypertrophic scars are not immediately visible after injury. These scars develop in 1 to 2 months after injury, whereas keloids develop months to years after the initial injury, supporting our results.
Experimental evidence suggests that fibroblasts from hypertrophic scars might represent a hyper-proliferative phenotype that can be reverted once the stimulation, such as the overabundance of growth factors, is lifted \cite{Tuan1998}. We verified this by setting the signaling molecule density to equilibrium on day 1000, and saw that this directly initiates the change of the (myo)fibroblast distributions and collagen density toward the healthy equilibria. First, the myofibroblasts disappear after 100 days, then 350 days later, the collagen density reaches equilibrium, and finally 100 days after that, the fibroblast distribution reaches equilibrium. Hence, according to our simulation, it takes about 1.5 year to reverse the process. To conclude this Section, the model is stable under the condition that the decay rate of the signaling molecules is not too far decreased to values below the stated bound $\delta_c\geq k_c/(a_c^{I}\overline{\rho})$. 

\section{Conclusion and discussion}
\label{sec:6}
In this study, we investigate the stability of the one-dimensional model for intensity of contraction and the formation of contractures in burn scars. The model presented in this paper is the one-dimensional version of the morphoelastic model developed by Koppenol. This model is based on the theory and derivations developed by \cite{Hall}. In this model, four constituents are incorporated: fibroblasts, myofibroblasts, signaling molecules, and collagen. Furthermore, we use equations for the displacement of the dermal layer, the displacement velocity of the dermal layer, and the effective Eulerian strain present in the dermal layer.

We presented a stability analysis for the model, for both the fully continuous and the semi-discrete (where the spatial derivatives have been replaced with differences) version of the problem. A surprising result was that we could derive the eigenvalues of the matrix involved in the stability analytically. This is possible because the linearised equations \eqref{Linearised} leave out other variables after accounting for the equilibria values. As a result, we could say that three of the six eigenvalues meet the stability constraints independent of the chosen value for the parameters, given that the parameters involved are positive and realistic. We have shown that the equilibrium distribution of the effective strain should should meet $\overline{\varepsilon}\leq 1$, and that the parameter that represents the viscosity of skin should be greater or equal to a factor containing the total mass density of dermal tissues, the Young’s Modulus and the equilibrium distribution of the effective strain, to have monotonically behavior of the solution. Note that the the stability criterium of the effective strain is also a physical requirement, from equation \eqref{PDE_e}. Further, another important stability constraint states the model is stable on the condition that the decay rate of the signaling molecules is greater than a factor concerning the maximum net secretion rate of the signaling molecules, the concentration of the signaling molecules that causes half-maximum net secretion rate of the signaling molecules, and the collagen equilibrium density.

We have shown that there is a consistency between the eigenvalues of the discrete model which we used for the uniform grid based finite element approximation, and the eigenvalues of the continuous model which is the `true’ model. We see that if the equilibrium solution to the continuous problem is stable, then the equilibrium to the semi-discrete problem is also stable under the current discretization (that is, if we use the right discretization method). The convergence rate towards the equilibrium is determined by the obtained eigenvalues of the system. 
In case $\mu \geq \frac{h}{\tan(\pi\beta h)} \sqrt{\rho_tE\sqrt{\overline{\rho}}(1-\overline{\varepsilon})}$, convergence in the semi-discrete system is monotonic for $\beta=1$ and consistent for $h\to0$. For monotonic convergence in the continuous system it must hold that $\mu\geq\frac{1}{\pi}\sqrt{\rho_tE\sqrt{\overline{\rho}}(1-\overline{\varepsilon})}\geq \frac{h}{\tan(\pi\beta h)} \sqrt{\rho_tE\sqrt{\overline{\rho}}(1-\overline{\varepsilon})}$, $\beta = 1,\ldots,n-1$. Hence monotonic convergence in the continuous system implies monotonic convergence in the semi-discrete system. Conversely, convergence could be monotonic in the semi-discrete system and not in the continuous system.
We have assessed the convergence of the numerical method experimentally, in which the order of convergence is shown to be of order $\mathcal{O}(h^2)$. Since the difference between the eigenvalues from the continuous and semi-discrete problem is of the order $\mathcal{O}(h^2)$, the convergence rates towards the equilibrium differ by an order $\mathcal{O}(h^2)$. This is in accordance to the expectations since the discretization method should have local truncation errors of order $\mathcal{O}(h^2)$. 

Using numerical simulations, we validated the stability constraints that we derived in the analysis. In case we meet stability criteria, the model behaves absolutely stable given these stable parameter values. Because of the initial perturbation, it takes some time to rearrange the distribution of (myo)fibroblasts and the densities of signaling molecules and collagen. First, the signaling molecule density increases in the center of the modeled skin, after which local fibroblasts differentiate to myofibroblasts, decreasing the local fibroblast distribution and increasing the local myofibroblast distribution. Because of the initial perturbation in the displacement velocity density, there are local displacements. Because the displacement velocity density rearranges such that all values have the same sign, the mesh moves in one particular direction. Both the collagen density and effective strain density gradually move toward the equilibrium. We conclude that a small perturbation of order $\mathcal{O}(10^{-15})$ g/cm$^3$ in the signaling molecule density and a few cells in the (myo)fibroblast distributions is already responsible for initializing wound healing that takes more than a year time.

In case we do not meet the signaling molecule stability condition $\delta_c\geq k_c/(a_c^{I}\overline{\rho})$, the model can numerically be unstable. Initially, the model seems stable. The signaling molecule density and myofibroblast distribution seem to reach equilibria first, after which the fibroblast distribution and the collagen density seem to reach equilibria as well. Shortly after this has happened, the signaling molecule density or the myofibroblast distribution increase first, after which the fibroblast distribution drops and the collagen density explodes. Although the signaling molecule density and (myo)fibroblast distributions move back towards the equilibria, the collagen density does not, and therefore the numerical method does not converge. Because of the lack of convergence in the inner Picard iterations, the numerical method cannot attain a solution.

We confirmed the model is stable if the eigenvalues are not real-valued. If the viscosity is low, the figures show that the distributions and densities of the variables reach their equilibrium densities. Though, in this case convergence is not monotonic, but oscillates as seen in Figure \ref{Stable2}. Besides this conclusion, we point out that the larger the number of initial perturbed waves, the faster the equilibria are reached and the faster the initial oscillations fade out. Because of an initial increase in amplitude in, and the oscillating behavior of, the displacement velocity density, the mesh moves shortly after the start of the simulation back and forth to the left and right. After the displacement velocity density stabilises, the distributions and densities of the constituents move gradually toward the equilibria. In conclusion, we need real-valued eigenvalues to prevent the model to increase the amplitudes of the initial perturbations in the displacement velocity density. However, this does not induce instability in terms of equilibria. 

If we have $\delta_c<k_c/(a_c^{I}\overline{\rho})$ not too far below the bound, then the signaling molecules move away from equilibrium and affect the distributions of the fibroblasts and the myofibroblasts. All the constituents move away from the expected equilibria and oscillate around new equilibria. The collagen density still shows the initial waves of the perturbations around day 260, while these waves already vanished in the other densities. We have linked this situation to real-life scar occurrences, namely hypertrophic scars and keloids. By reverting the stimulation of matrix production and differentiation to myofibroblasts by setting the signaling molecule density to (healthy) equilibrium, we have provided experimental evidence, from a mathematical point of view, that fibroblasts can be indeed be reverted. 
Taken together, the numerical model fully reproduces the stability constraints.

It would interest to incorporate hypertrophy to this one-dimensional morphoelastic model \cite{Koppenol2017b}, since hypertrophic scars can also develop contractures and elevate above healthy skin levels. The beauty of the one-dimensional model is the speed, hence incorporation of hypertrophy will quickly yield new results and therefore insight. However, validating results from such a model is a challenge since hypertrophy depends highly on angiogenesis, which nowadays seems impossible to test in vitro.

An interesting direction is to model the boundaries of the wounded area as elastic springs, since with the current setting the boundary of the domain of computation needs to be ``sufficiently far away". We are planning on incorporating pulling and stretching forces because of the growth of children and motility. A first attempt to incorporate the growth of children is to add another term to the right-hand side of equation \eqref{PDE_e}. We will incorporate motility forces by adding new boundary conditions.

Considering the modeling choices, we could keep a linear growth rate and introducing a tune-able quadratic cell death term for fitting equilibrium, instead of with the constant $q$ in equation \eqref{PDE_N}. We can also easily consider that myofibroblasts in response to TGF-$\beta$ move slower than fibroblasts \cite{Thampatty2006}, and that myofibroblasts can differentiate back to fibroblasts under the influence of Prostaglandin E2 (PGE2) \cite{Garrison2013}.

It would interest to investigate the influence of such other modeling choices on the simulation results.
Currently, we are working on a sensitivity analysis and feasibility study, and a neural network for the one-dimensional morphoelastic model.

\section*{Acknowledgement}
The authors are grateful for the financial support by the Dutch Burns Foundation under Project 17.105.

\section*{Conflict of interest}
The authors declare that they have no conflict of interest.

\bibliographystyle{apa}
\bibliography{references}

\section*{Appendices}
\label{sec:7}
\subsection*{Appendix 1: The derivation of the stability constraints for the continuous problem}
First we substitute the variations \eqref{Variations_continuous} into the linearised equations \eqref{Linearised}. This yields
\begin{multline*}
\frac{1}{|\Omega|}\sum_{j = -\infty}^{\infty} \dot{c}^c_j(t) e^{2 i \pi j x}
+ \frac{D_c}{|\Omega|}\sum_{j = -\infty}^{\infty} (2\pi j)^2 c^c_j(t) e^{2 i \pi j x} \\
+\frac{\overline{N}}{|\Omega|}\left[  \delta_c\overline{\rho}- \frac{ k_c }{ a_c^{I} } \right]\sum_{j = -\infty}^{\infty} c^c_j(t) e^{2 i \pi j x}=0,
\end{multline*}
\begin{multline*}\frac{1}{|\Omega|}\sum_{j = -\infty}^{\infty} \dot{c}^N_j(t) e^{2 i \pi j x}
+ \frac{D_F\overline{N}}{|\Omega|}\sum_{j = -\infty}^{\infty} (2\pi j)^2 c^N_j(t) e^{2 i \pi j x}\\
- \frac{\chi_F\overline{N}}{|\Omega|}\sum_{j = -\infty}^{\infty} (2\pi j)^2 c^c_j(t) e^{2 i \pi j x}
+\frac{\delta_N}{|\Omega|}\sum_{j = -\infty}^{\infty} c^N_j(t) e^{2 i \pi j x}\\
- \frac{r_F}{|\Omega|}\overline{N}^q((1+q)(1-\kappa_F\overline{N})-\kappa_F\overline{N})\sum_{j = -\infty}^{\infty} c^N_j(t) e^{2 i \pi j x} \\
+\frac{r_F\kappa_F\overline{N}^{1+q}}{|\Omega|}\sum_{j = -\infty}^{\infty} c^M_j(t) e^{2 i \pi j x}\\
- \frac{\overline{N}}{|\Omega|}\left[\frac{r_Fr_F^{\text{max}}}{a_c^{III}}[1-\kappa_F\overline{N}]\overline{N}^q-k_F\right] \sum_{j = -\infty}^{\infty} c^c_j(t) e^{2 i \pi j x}=0,
\end{multline*}
\begin{multline*}\frac{1}{|\Omega|}\sum_{j = -\infty}^{\infty} \dot{c}^M_j(t) e^{2 i \pi j x}
+ \frac{D_F\overline{N}}{|\Omega|}\sum_{j = -\infty}^{\infty} (2\pi j)^2 c^M_j(t) e^{2 i \pi j x} \\
+ \frac{\delta_M}{|\Omega|}\sum_{j = -\infty}^{\infty} c^M_j(t) e^{2 i \pi j x} 
- \frac{k_F\overline{N}}{|\Omega|}\sum_{j = -\infty}^{\infty} c^c_j(t) e^{2 i \pi j x} = 0,
\end{multline*}
\begin{multline*}\frac{1}{|\Omega|}\sum_{j = -\infty}^{\infty} \dot{c}^{\rho}_j(t) e^{2 i \pi j x}
+\frac{\delta_\rho\overline{\rho}^2}{|\Omega|}(\eta^{II}-\eta^I)\sum_{j = -\infty}^{\infty} c^M_j(t) e^{2 i \pi j x}\\
-\frac{\delta_\rho\overline{\rho}^2\overline{N}}{|\Omega|}\left(\frac{k_\rho^{max}}{a_c^{IV}}+a_c^{II}\right)\sum_{j = -\infty}^{\infty} c^c_j(t) e^{2 i \pi j x}\\
+\frac{2\delta_\rho\overline{N}}{|\Omega|}\overline{\rho}\sum_{j = -\infty}^{\infty} c^{\rho}_j(t) e^{2 i \pi j x}=0,
\end{multline*}
for the chemical part of the model, and 
\begin{multline*}\frac{1}{|\Omega|}\sum_{j = -\infty}^{\infty} \dot{c}^v_j(t) e^{2 i \pi j x}
+ \frac{\mu}{|\Omega|\rho_t}\sum_{j = -\infty}^{\infty} (2\pi j)^2 c^v_j(t) e^{2 i \pi j x} \\
- i\frac{E \sqrt{\overline{\rho}}}{|\Omega|\rho_t}\sum_{j = -\infty}^{\infty} (2\pi j) c^{\varepsilon}_j(t) e^{2 i \pi j x}
- i\frac{E\overline{\varepsilon}}{|\Omega|2\rho_t\sqrt{\overline{\rho}}}\sum_{j = -\infty}^{\infty} (2\pi j) c^{\rho}_j(t) e^{2 i \pi j x}\\
- i\frac{\xi\overline{\rho}}{|\Omega|\rho_t(R^2+\overline{\rho}^2)}\sum_{j = -\infty}^{\infty} (2\pi j) c^M_j(t) e^{2 i \pi j x}=0,
\end{multline*}
\begin{multline*}
\frac{1}{|\Omega|}\sum_{j = -\infty}^{\infty} \dot{c}^{\varepsilon}_j(t) e^{2 i \pi j x}
+ i\frac{(\overline{\varepsilon}-1)}{|\Omega|}\sum_{j = -\infty}^{\infty} (2\pi j) c^v_j(t) e^{2 i \pi j x}\\
+ \frac{\zeta\overline{\varepsilon}\overline{N}}{|\Omega|}\sum_{j = -\infty}^{\infty} c^c_j(t) e^{2 i \pi j x}=0,
\end{multline*}
for the mechanical part of the model. Multiplication by $e^{-2i\pi kx}$ gives
\begin{multline*}\frac{1}{|\Omega|}\sum_{j = -\infty}^{\infty} \dot{c}^c_j(t) e^{2 i \pi (j-k) x}
+ \frac{D_c}{|\Omega|}\sum_{j = -\infty}^{\infty} (2\pi j)^2 c^c_j(t) e^{2 i \pi (j-k) x} \\
+\frac{\overline{N}}{|\Omega|}\left[  \delta_c\overline{\rho}- \frac{ k_c }{ a_c^{I} } \right]\sum_{j = -\infty}^{\infty} c^c_j(t) e^{2 i \pi (j-k) x}=0,
\end{multline*}
\begin{multline*}
\frac{1}{|\Omega|}\sum_{j = -\infty}^{\infty} \dot{c}^N_j(t) e^{2 i \pi (j-k) x}
+ \frac{D_F\overline{N}}{|\Omega|}\sum_{j = -\infty}^{\infty} (2\pi j)^2 c^N_j(t) e^{2 i \pi (j-k) x}\\
- \frac{\chi_F\overline{N}}{|\Omega|}\sum_{j = -\infty}^{\infty} (2\pi j)^2 c^c_j(t) e^{2 i \pi (j-k) x}\\
- \frac{r_F}{|\Omega|}\overline{N}^q((1+q)(1-\kappa_F\overline{N})-\kappa_F\overline{N})\sum_{j = -\infty}^{\infty} c^N_j(t) e^{2 i \pi (j-k) x} \\
+\frac{\delta_N}{|\Omega|}\sum_{j = -\infty}^{\infty} c^N_j(t) e^{2 i \pi (j-k) x}
+\frac{r_F\kappa_F\overline{N}^{1+q}}{|\Omega|}\sum_{j = -\infty}^{\infty} c^M_j(t) e^{2 i \pi (j-k) x}\\
- \frac{\overline{N}}{|\Omega|}\left[\frac{r_Fr_F^{\text{max}}}{a_c^{III}}[1-\kappa_F\overline{N}]\overline{N}^q-k_F\right] \sum_{j = -\infty}^{\infty} c^c_j(t) e^{2 i \pi (j-k) x}=0,
\end{multline*}
\begin{multline*}
\frac{1}{|\Omega|}\sum_{j = -\infty}^{\infty} \dot{c}^M_j(t) e^{2 i \pi (j-k) x}
+ \frac{D_F\overline{N}}{|\Omega|}\sum_{j = -\infty}^{\infty} (2\pi j)^2 c^M_j(t) e^{2 i \pi (j-k) x} \\
+ \frac{\delta_M}{|\Omega|}\sum_{j = -\infty}^{\infty} c^M_j(t) e^{2 i \pi (j-k) x} 
- \frac{k_F\overline{N}}{|\Omega|}\sum_{j = -\infty}^{\infty} c^c_j(t) e^{2 i \pi (j-k) x} = 0,
\end{multline*}
\begin{multline*}
\frac{1}{|\Omega|}\sum_{j = -\infty}^{\infty} \dot{c}^{\rho}_j(t) e^{2 i \pi (j-k) x}
+\frac{\delta_\rho\overline{\rho}^2}{|\Omega|}(\eta^{II}-\eta^I)\sum_{j = -\infty}^{\infty} c^M_j(t) e^{2 i \pi (j-k) x}\\
-\frac{\delta_\rho\overline{\rho}^2\overline{N}}{|\Omega|}\left(\frac{k_\rho^{max}}{a_c^{IV}}+a_c^{II}\right)\sum_{j = -\infty}^{\infty} c^c_j(t) e^{2 i \pi (j-k) x}\\
+\frac{2\delta_\rho\overline{N}}{|\Omega|}\overline{\rho}\sum_{j = -\infty}^{\infty} c^{\rho}_j(t) e^{2 i \pi (j-k) x}=0,
\end{multline*}
for the chemical part of the model, and
\begin{multline*}
\frac{1}{|\Omega|}\sum_{j = -\infty}^{\infty} \dot{c}^v_j(t) e^{2 i \pi (j-k) x}
+ \frac{\mu}{|\Omega|\rho_t}\sum_{j = -\infty}^{\infty} (2\pi j)^2 c^v_j(t) e^{2 i \pi (j-k) x} \\
- i\frac{E \sqrt{\overline{\rho}}}{|\Omega|\rho_t}\sum_{j = -\infty}^{\infty} (2\pi j) c^{\varepsilon}_j(t) e^{2 i \pi (j-k) x}
- i\frac{E\overline{\varepsilon}}{|\Omega|2\rho_t\sqrt{\overline{\rho}}}\sum_{j = -\infty}^{\infty} (2\pi j) c^{\rho}_j(t) e^{2 i \pi (j-k) x}\\
- i\frac{\xi\overline{\rho}}{|\Omega|\rho_t(R^2+\overline{\rho}^2)}\sum_{j = -\infty}^{\infty} (2\pi j) c^M_j(t) e^{2 i \pi (j-k) x}=0,
\end{multline*}
\begin{multline*}
\frac{1}{|\Omega|}\sum_{j = -\infty}^{\infty} \dot{c}^{\varepsilon}_j(t) e^{2 i \pi (j-k) x}
+ i\frac{(\overline{\varepsilon}-1)}{|\Omega|}\sum_{j = -\infty}^{\infty} (2\pi j) c^v_j(t) e^{2 i \pi (j-k) x}\\
+ \frac{\zeta\overline{\varepsilon}\overline{N}}{|\Omega|}\sum_{j = -\infty}^{\infty} c^c_j(t) e^{2 i \pi (j-k) x}=0,
\end{multline*}
for the mechanical part of the model. Integration over $\Omega$ gives the result, hence equations \eqref{Integration_continuous1} and \eqref{Integration_continuous2}.

\subsection*{Appendix 2: The derivation of the stability constraints for the discrete problem}
Substitution of the variations \eqref{Variations_discrete} in finite differences equations \eqref{lambda1} and \eqref{lambda2} gives
\begin{multline*}
\lambda c_k =
- \frac{D_c}{h^2}\sum_{\beta=1}^{n-1}\hat{c}_\beta\left\{e^{-2\pi\beta (k-1)hi} -2 e^{-2\pi\beta khi} + e^{-2\pi\beta (k+1)hi}\right\} \\
+\overline{N}\left[  \delta_c\overline{\rho}- \frac{ k_c }{ a_c^{I} } \right]\sum_{\beta=1}^{n-1} \hat{c}_\beta e^{-2\pi\beta khi},
\end{multline*}
\begin{multline*}
\lambda N_k = - \frac{D_F\overline{N}}{h^2}\sum_{\beta=1}^{n-1}\hat{N}_\beta\left\{e^{-2\pi\beta (k-1)hi} - 2 e^{-2\pi\beta khi} + e^{-2\pi\beta (k+1)hi}\right\} \\
+ \frac{\chi_F\overline{N}}{h^2}\sum_{\beta=1}^{n-1}\hat{c}_\beta\left\{e^{-2\pi\beta (k-1)hi} -2 e^{-2\pi\beta khi} + e^{-2\pi\beta (k+1)hi}\right\}\\
+\left[\delta_N - r_F\overline{N}^q((1+q)(1-\kappa_F\overline{N})-\kappa_F\overline{N})\right]\sum_{\beta=1}^{n-1} \hat{N}_\beta e^{-2\pi\beta khi}\\
+r_F\kappa_F\overline{N}^{1+q}\sum_{\beta=1}^{n-1} \hat{M}_\beta e^{-2\pi\beta khi}\\
- \overline{N}\left[\frac{r_Fr_F^{\text{max}}}{a_c^{III}}[1-\kappa_F\overline{N}]\overline{N}^q-k_F\right]\sum_{\beta=1}^{n-1} \hat{c}_\beta e^{-2\pi\beta khi}, 
\end{multline*}
\begin{multline*}
\lambda M_k = 
- \frac{D_F\overline{N}}{h^2}\sum_{\beta=1}^{n-1}\hat{M}_\beta\left\{e^{-2\pi\beta (k-1)hi} -2 e^{-2\pi\beta khi} + e^{-2\pi\beta (k+1)hi}\right\}\\
+\delta_M\sum_{\beta=1}^{n-1} \hat{M}_\beta e^{-2\pi\beta khi} - k_F\overline{N}\sum_{\beta=1}^{n-1} \hat{c}_\beta e^{-2\pi\beta khi},
\end{multline*}
\begin{multline*}
\lambda \rho_k =
\delta_\rho\overline{\rho}^2(\eta^{II}-\eta^I)\sum_{\beta=1}^{n-1} \hat{M}_\beta e^{-2\pi\beta khi} \\
- \delta_\rho\overline{\rho}^2\overline{N}\left(\frac{k_\rho^{max}}{a_c^{IV}}+a_c^{II}\right)\sum_{\beta=1}^{n-1} \hat{c}_\beta e^{-2\pi\beta khi} 
+ 2\delta_\rho\overline{N}\overline{\rho}\sum_{\beta=1}^{n-1} \hat{\rho}_\beta e^{-2\pi\beta khi},
\end{multline*}
for the chemical part of the model, and
\begin{multline*}
\lambda v_k= -\frac{\mu}{\rho_t h^2}\sum_{\beta=1}^{n-1}\hat{v}_\beta\left\{e^{-2\pi\beta (k-1)hi} - 2 e^{-2\pi\beta khi} + e^{-2\pi\beta (k+1)hi}\right\} \\
- \frac{E\sqrt{\overline{\rho}}}{\rho_t2h}\sum_{\beta=1}^{n-1}\hat{\varepsilon}_\beta\left\{ e^{-2\pi\beta (k+1)hi} -e^{-2\pi\beta (k-1)hi}\right\}\\
- \frac{E\overline{\varepsilon}}{2\rho_t\sqrt{\overline{\rho}}2h}\sum_{\beta=1}^{n-1}\hat{\rho}_\beta\left\{ e^{-2\pi\beta (k+1)hi} - e^{-2\pi\beta (k-1)hi}\right\}\\
-\frac{\xi\overline{\rho}}{\rho_t(R^2+\overline{\rho}^2)2h}\sum_{\beta=1}^{n-1}\hat{M}_\beta\left\{  e^{-2\pi\beta (k+1)hi}-  e^{-2\pi\beta (k-1)hi}\right\},
\end{multline*}
\begin{equation*}
\lambda\varepsilon_k = \frac{(\overline{\varepsilon}-1)}{2h}\sum_{\beta=1}^{n-1}\hat{v}_\beta\left\{e^{-2\pi\beta (k+1)hi}- e^{-2\pi\beta (k-1)hi}\right\}
+\zeta\overline{\varepsilon}\overline{N}\sum_{\beta=1}^{n-1} \hat{c}_\beta e^{-2\pi\beta khi},
\end{equation*}
for the mechanical part of the model. 

This must be true for arbitrary $\{c_\beta,N_\beta,M_\beta,\rho_\beta,v_\beta,\varepsilon_\beta\}$, hence each factor following $\{c_\beta,N_\beta,M_\beta,\rho_\beta,v_\beta,\varepsilon_\beta\}$ in the sum should be zero.
Subdivision by $e^{-2\pi\beta khi}$ gives
\begin{equation*}
\lambda c_k =
- \frac{D_c}{h^2}\hat{c}_\beta\left\{e^{2\pi\beta hi} -2 + e^{-2\pi\beta hi}\right\} 
+\overline{N}\left[  \delta_c\overline{\rho}- \frac{ k_c }{ a_c^{I} } \right]\hat{c}_\beta,
\end{equation*}
\begin{multline*}
\lambda N_k =  
- \frac{D_F\overline{N}}{h^2}\hat{N}_\beta\left\{e^{2\pi\beta hi} - 2 + e^{-2\pi\beta hi}\right\} \\
+ \frac{\chi_F\overline{N}}{h^2}\hat{c}_\beta\left\{e^{2\pi\beta hi} -2 + e^{-2\pi\beta hi}\right\}\\
+\left[\delta_N - r_F\overline{N}^q((1+q)(1-\kappa_F\overline{N})-\kappa_F\overline{N})\right]\hat{N}_\beta\\
+r_F\kappa_F\overline{N}^{1+q}\hat{M}_\beta
- \overline{N}\left[\frac{r_Fr_F^{\text{max}}}{a_c^{III}}[1-\kappa_F\overline{N}]\overline{N}^q-k_F\right]\hat{c}_\beta, 
\end{multline*}
\begin{equation*}
\lambda M_k= 
- \frac{D_F\overline{N}}{h^2}\hat{M}_\beta\left\{e^{2\pi\beta hi} -2 + e^{-2\pi\beta hi}\right\}
+\delta_M\sum_{\beta=1}^{n-1} \hat{M}_\beta - k_F\overline{N}\hat{c}_\beta,
\end{equation*}
\begin{equation*}
\lambda \rho_k=
\delta_\rho\overline{\rho}^2(\eta^{II}-\eta^I)\hat{M}_\beta
- \delta_\rho\overline{\rho}^2\overline{N}\left(\frac{k_\rho^{max}}{a_c^{IV}}+a_c^{II}\right)\hat{c}_\beta
+ 2\delta_\rho\overline{N}\overline{\rho}\hat{\rho}_\beta,
\end{equation*}
for the chemical part of the model, and
\begin{multline*}
\lambda v_k = -\frac{\mu}{\rho_th^2}\hat{v}_\beta\left\{e^{2\pi\beta hi} - 2 + e^{-2\pi\beta hi}\right\} 
- \frac{E\sqrt{\overline{\rho}}}{\rho_t2h}\hat{\varepsilon}_\beta\left\{ e^{-2\pi\beta hi} -e^{2\pi\beta hi}\right\}\\
- \frac{E\overline{\varepsilon}}{2\rho_t\sqrt{\overline{\rho}}2h}\hat{\rho}_\beta\left\{ e^{-2\pi\beta hi} - e^{2\pi\beta hi}\right\}\\
-\frac{\xi\overline{\rho}}{\rho_t(R^2+\overline{\rho}^2)2h}\hat{M}_\beta\left\{  e^{-2\pi\beta hi}-  e^{2\pi\beta hi}\right\},
\end{multline*}
\begin{equation*}
\lambda\varepsilon_k= \frac{(\overline{\varepsilon}-1)}{2h}\hat{v}_\beta\left\{e^{-2\pi\beta hi}- e^{2\pi\beta hi}\right\}
+\zeta\overline{\varepsilon}\overline{N}\hat{c}_\beta,
\end{equation*}
for the mechanical part of the model. Using Euler's formula and $2-2\cos(2\pi\beta h)=4\sin^2(\pi\beta h)$ gives the result, hence equations \eqref{Substitution_discrete1} and \eqref{Substitution_discrete2}.

\subsection*{Appendix 3: The parameters}
The value of the parameter $q$ is a consequence of the values of other parameters, see equation \eqref{delta_N}. 
Further, the value of the parameter $k_\rho$ is a consequence of the values of other parameters, see equation \eqref{rho_eq}.


\begin{table}
\caption{Overview of the parameters used for the simulations. Shown are the symbols, values, dimensions and references. Here TW denotes that the value of the parameter is estimated in the study, and NC denotes that the value of the parameter is a consequence because of the chosen values for other parameters.}
\label{Parameters}
\begin{tabular}{llll} 
\hline\noalign{\smallskip}
Symbol 				& Value & Dimension & Reference\\
\noalign{\smallskip}\hline\noalign{\smallskip}
$D_c$ 				& $2.88\times10^{-3}$ 	& cm$^2$/day 				&\cite{Haugh2006}\\
$D_F$ 				& $10^{-7}$ 				& cm$^5$/(cells day) 			&\cite{Sillman2003}\\
$\chi_F$ 				& $2\times10^{-3}$ 		& cm$^5$/(g day) 			&\cite{Murphy2012}\\
$k_c$ 				& $4\times10^{-13}$ 		& g/(cells day) 				&\cite{Olsen1995}\\
$r_F$ 				& $9.24\times10^{-1}$ 	& cm$^{3q}$/(cells$^q$ day) 	&\cite{Alberts1989} \&\cite{Gosh}\\
$r_F^{\text{max}}$ 	& $2$ 					& - 							&\cite{Strutz2001}\\
$k_\rho$ 				& $7.6\times10^{-8}$ 		& g/(cells day) 				&\text{[NC]}\\
$k_\rho^{\text{max}}$ & $10$ 					& - 							&\cite{Olsen1995}\\
$a_c^{I}$ 			& $10^{-8}$ 				& g/cm$^3$ 					&\cite{Olsen1995}\\
$a_c^{II}$ 			& $2\times10^8$ 			& cm$^3$/g 					&\cite{Overall1991}\\
$a_c^{III}$ 			& $10^{-8}$ 				& g/cm$^3$ 					&\cite{Grotendorst1992} \&\cite{Olsen1995}\\
$a_c^{IV}$ 			& $10^{-9}$ 				& g/cm$^3$ 					&\cite{Roberts1986}\\
$\eta^I$ 			& $2$ 					& - 							&\cite{Rudolph1991}\\
$\eta^{II}$ 			& $5\times10^{-1}$ 		& - 							&\text{[TW]}\\
$k_F$ 				& $1.08\times10^{7}$ 		& cm$^3$/(g day) 			&\cite{Desmouliere1993}\\
$\kappa_F$ 			& $10^{-6}$ 				& cm$^3$/cells 				&\cite{VandeBerg1989}\\
$q$ 					& $-4.151\times10^{-1}$ 	& -  							&\text{[NC]}\\
$\delta_c$ 			& $5\times10^{-4}$ 		& cm$^6$/(cells g day)		 	&\cite{Olsen1995}\\
$\delta_N$ 			& $2\times10^{-2}$ 		& /day 						&\cite{Olsen1995}\\
$\delta_M$ 			& $6\times10^{-2}$ 		& /day 						&\cite{Koppenol2017b}\\
$\delta_\rho$ 			& $6\times10^{-6}$ 		& cm$^6$/(cells g day) 		&\cite{Koppenol2017b}\\
$\overline{N}$ 		& $10^{4}$ 				& cells/cm$^3$ 				&\cite{Olsen1995}\\
$\overline{M}$ 		& $0$ 					& cells/cm$^3$ 				&\cite{Olsen1995}\\
$\overline{c}$ 		& $0$ 					& g/cm$^3$ 					&\cite{Koppenol2017b}\\
$\overline{\rho}$ 		& $1.125\times10^{-1}$ 	& g/cm$^3$ 					&\cite{Olsen1995}\\
$\rho_t$ 				& 1.09 					& g/cm$^3$ 					&\cite{Wrobel2009}\\
$\mu$ 				& $10^2$ 				& (N day)/cm$^2$  			&\text{[TW]}\\
$E$ 					& $2.1\times10^2$ 		& N/((g cm)$^{0.5})$ 			&\text{[TW]}\\
$\xi$ 				& $4.4\times10^{-2}$ 		& (N g)/(cells cm$^2$) 			& \cite{Maskarinec} \& \cite{Wrobel}\\
$R$ 					& $9.95\times10^{-1}$ 	& g/cm$^3$ 					& \text{[TW]}\\
$\zeta$ 				& $4\times10^2$ 			& cm$^6$/(cells g day) 		&\text{[TW]} \\
\noalign{\smallskip}\hline
\end{tabular}
\end{table}

\subsection*{Appendix 4: Absolute errors in convergence}
The averaged errors in Table \ref{tab:slopes} show that the order of convergence in the numerical method is $\mathcal{O}(h^2)$. Shown are the averaged slopes of the errors that are defined in \eqref{norm_e41} and by:
\begin{align*}
\epsilon_{L^1}(h) &= h\sum_{i=1}^n     \left|z_{h/r}(x_{i,n}) - z_h(x_{i,n})\right|,\\
\epsilon_{L^2}(h) &= \sqrt{h\sum_{i=1}^n     \left(z_{h/r}(x_{i,n}) - z_h(x_{i,n})\right)^2},
\end{align*}
where the grid-points $x_{i,n}$ correspond to the grid-points in the simulation with $n$ nodes.
\begin{table}
\caption{Overview of the averaged slopes of the errors of the variables for different element sizes $h$ on the full domain of computation and on the boundary of the wound. The columns show slopes for the different errors, and the rows show the averaged slopes for the variables. The last column shows the averaged slopes of the rows. The reference is the solution in which $h=0.0078$.}
\label{tab:slopes}
\begin{tabular}{l|llll|l} 
\hline\noalign{\smallskip}
Variable 		& $\epsilon_{|41|}$ & $\epsilon_{L^1}$ & $\epsilon_{L^2}$ & $\epsilon_{boundary}$ & Averaged\\
\noalign{\smallskip}\hline\noalign{\smallskip}
$N$			&2.1843 & 2.0160 & 1.9701 & 2.1850 & 2.0889\\ 
$M$			&2.1735 & 2.1203 & 2.0961 & 2.1892 & 2.1448\\ 
$c$			&2.1900 & 2.0964 & 2.0675 & 2.0929 & 2.1117\\ 
$\rho$		&2.1911 & 2.0626 & 1.9211 & 2.1708 & 2.0864\\ 
$v$			&2.1882 & 2.1891 & 2.1911 & 2.1189 & 2.1718\\ 
$\varepsilon$	&2.2283 & 2.2301 & 2.2521 & 2.2403 & 2.2377\\
\noalign{\smallskip}\hline
\end{tabular}
\end{table}
\end{document}